\documentclass[12pt]{amsart}
\usepackage{graphicx}
\usepackage{amssymb,latexsym,amsmath}
\newcommand {\PP}{{I\kern-.3em P}}
\newcommand {\ZZ}{{Z\kern-.45em Z}}
\newcommand {\RR}{{\mathbb{R}}}
\newcommand {\HH}{\Bbb{H}}
\newcommand {\CC}{{\mathbb{C}}}
\newcommand {\NN}{\Bbb{N}}

\newcommand{\beq}{\begin{equation}}
\newcommand{\eeq}{\end{equation}}

\newtheorem{thm}{Theorem}
\newtheorem{lemma}{Lemma}
\newtheorem{prop}{Proposition}
\newtheorem{cor}{Corollary}
\theoremstyle{remark}

\theoremstyle{definition}

\newcommand{\aal}{\alpha }
\def\d{\displaystyle}

\begin{document}

\title{The geometry at infinity of a hyperbolic Riemann surface of infinite type }
\author{Andrew Haas and Perry Susskind}

\begin{abstract}
We study geodesics on planar Riemann surfaces of infinite type having a single infinite end. Of particular interest is the class of geodesics that go out the infinite end in a most efficient manner. We investigate properties of these geodesics and relate them to the structure of the boundary of a   Dirichlet polygon for a Fuchsian group representing the surface.

\end{abstract}
\subjclass{ 30F35, 30F45, 53C22}
\date{}
\keywords{Geodesic, hyperbolic surface, Dirichlet polygon }

\maketitle
\markboth{Surfaces of infinite type}
{Surfaces of infinite type}
\section{introduction}

A flute surface $S$ is  most simply 
described as a connected domain in the complex plane, for which all but one of the components in the complement of $S$  
 is isolated from the others. Flute surfaces were first considered by Basmajian \cite{bas1} as examples of  the simplest sort of hyperbolic Riemann
surface of infinite type. A flute surface has a single infinite end.
  The presence of
such an infinite end, even one of this simple sort, allows for
many different possibilities for the geometry of the surface, which
have no parallels in the theory of finite surfaces.
 In this paper,   our main
concern is with the  behavior of the geometry associated to the infinite end
of the surface, as described by the special classes of infinite critical
and subcritical geodesic rays. These are geodesic rays that head
either directly, or almost directly, out the infinite end of the
surface. In the  theory of Fuchsian groups, these types of
rays are related to the existence of Dirichlet and Garnett points
in the limit set of a  Fuchsian group representing the surface
 \cite{haas, Waterman, Sullivan}.

Our approach is to employ a sequence of cut and paste operations to
construct flute surfaces with complex end structure, where the
building blocks are the simple untwisted flutes surfaces studied in
\cite{bas1, haas}.  We refer to the surfaces constructed in this way as {\it quilted} flute surfaces.  
The geometry out the infinite end of a quilted flute surface is considerably more complex than the end geometry
of an untwisted flute. Nevertheless, we demonstrate that there are certain fundamental similarities.

The boundary at infinity of a Dirichlet polygon for a Fuchsian group may be regarded as one measure of 
the complexity of the end geometry of the surface represented by the group. Applying our construction of quilted surfaces we show that, up to a sparse set, one can exercise surprising control over the boundary at infinity of 
a Dirichlet polygon for a quilted surface group. To this end we prove

\begin{thm}\label{C}
Let $K$ be a compact subset of $\RR$. There is a Fuchsian group $G$ representing a quilted flute surface and a point $\tilde{p}\in\HH^2$ so that
the boundary at infinity of the Dirichlet polygon for $G$ centered at $\tilde{p}$ consists of the union of the set $K$ and a countable set of isolated parabolic fixed points of $G$.
\end{thm}

The paper is structured as follows. We begin in Section \ref{review}, with a review of some facts about flute surfaces. In Section \ref{quilts} we construct  quilted flute surfaces and derive some of their basic properties.  In Section \ref{number4} we investigate a special class of canonical critical rays that go out the infinite end of a quilted flute surface.   The complexity of the set of such rays is one measure of the complexity of the surface.     Before moving on to probe more results about quilted surfaces we must turn to the hyperbolic plane. In Section \ref{lemmas} we prove several lemmas, crucial to the proofs in later sections.   
In Section \ref{example2} we develop an intrinsic version of Theorem \ref{C} for quilted surfaces and show how this can be used to  to prove Theorem \ref{C}.
 In the last section we look at  the finer structure of the set of infinite critical and subcritical geodesics rays. The main result of Section \ref{close} is that there is an underlying geodesic scaffolding, heading out the infinite end, which all infinite critical and subcritical geodesic rays must approach asymptotically. This is a broad generalization of a similar result about untwisted flutes, that appeared in \cite{haas}.
In obtaining these results we employ   a number of lemmas. Several of these are results about plane hyperbolic geometry which are interesting in their own right.
 %review%%%%%%%%%%%%%%%%%%%%%%%%%%%%%%%%%%%%%%%%%%%

\section{Basic properties of flute surfaces}\label{review}
The main reference for the material in this section is \cite{haas}.

Define an {\em end} {\sf E} of a manifold $M$ as follows. Let $
K_1 \subset K_{2} \subset ...\subset M$ be a nested sequence of
compact subsets of $M$ so that $\bigcup^{\infty}_{i=1} K_{i} = M$.
An end {\sf E}  is a sequence of connected components $\{{\sf E}_i\}$
in the complement of $K_{i}$ so that ${\sf E}_{i+1} \subset {\sf
E}_i$. This definition can be made independent of the given
exhaustion  $\{ K_i\}$. A   ray $\sigma$ is said to {\em go
out the end} {\sf E}, if for each integer $i > 0$ all but a
compact segment of $\sigma$ belongs to ${\sf E}_i.$

Let  $S$ be a hyperbolic surface. An end  {\sf E} of $S$ is called
a {\em puncture} if there is a subset $D$ of $S$ which is
conformally equivalent to the punctured disc $\{ z \ | \
0<|z|<1\}$ and for $i$ large, $ {\sf E}_i\subset D$. Similarly, we
call  {\sf E} a {\em hole} if there is a subset $D$ of $S$ which
is conformally equivalent to an annulus $\{ z \ |\ 1<|z|<r\}$ for
some $r>1$ and for $i$ large, $ {\sf E}_i\subset D$. {\sf E} is a
{\em finite end} if it is either a puncture or  a hole; otherwise
it is an {\em infinite end}.

An end ${\sf E} = \{{\sf E}_i\}$ is said to be of
the  {\em second kind} if $S$ contains a half-plane $P$ and for all $i>
0$, ${\sf E}_i\bigcap P\not=\emptyset$. If an end is not of the
second kind then it is of the {\em first kind}. A puncture is of the first
kind and a hole is
of the second kind.
 
Let ${\bf S}$ denote the infinite cylinder ${\bf C}={\rm
S}^1\times(0,\infty),$ with the set of points $\{ (1,n)|\,n\in \Bbb{N}
\}$ deleted, and define the space
 $\mathcal{ F}$ of isometry classes of complete metrics of constant
curvature  -1, that is, hyperbolic metrics on the surface {\sf  S}.
Define an involution $r:
{\sf S}\rightarrow {\sf S}$ by $r(e^{i\theta},t)
=(e^{-i\theta},t)$. Let $\mathcal{ F}_0 \subset \mathcal{  F}$ be the set of
isometry classes in $\mathcal{  F}$ for which there exists a
representative surface on which $r$ is an isometry.   Henceforth we
shall treat elements of $\mathcal{  F}$ and $\mathcal{  F}_0$ as hyperbolic
surfaces and suppose, in the latter case, that  $r$ is an
isometry. A
surface  in $\mathcal{  F}$ is called a {\em flute}; one in $\mathcal{  F}_0$ is
called an {\em untwisted flute}. An explicit construction of   flutes is given in \cite{bas1}. A flute has one infinite end and has 
a finite end associated to each of the deleted points $(1,n)$ and
to the ideal boundary ${\rm S}^1\times\{0\}.$ Note that each of these ends can be of the first or  the second kind, depending on  the hyperbolic metric we have chosen.
The proposition below, whose proof is given after we develop some notation, shows that this definition of a flute surface is consistent with the
 definition given in the introduction.

\begin{prop}\label{planar}
$S$ is a flute surface if and only if $S$ is conformally equivalent to a connected domain in $\CC$ with a single infinite end.
\end{prop}
 Let $S$ be a hyperbolic surface and $\sigma :[0,\infty)\rightarrow
S$ a geodesic ray.  Here and henceforth, all geodesics   are
parameterized by arc length. Define the function
$\Delta_{\sigma}(t)= t - d_{S}(\sigma(0),\sigma(t)) $, where
$d_{S}$ denotes distance as defined by the hyperbolic metric on
$S$. The ray $\sigma$ is then said to be {\em horocyclic,
critical, {\em or} subcritical}\/ if $\Delta_{\sigma}$ is
respectively, unbounded, zero, or nonzero but bounded. A critical
ray may be said to  travel directly out an end of
$S,$ and a subcritical ray may be said to  travel almost
directly out an end. It is known (see, \cite{haas, Waterman}) that critical rays
are simple and subcritical rays are simple beyond some point.

Let ${\sf E} = \{{\sf E}_i\}$ be an end of a flute $F$. It is
known from \cite{haas} that  for any $p\in F$ there is a critical
ray with initial point $p$ that goes out the end {\sf E}. When the
end ${\sf E}$ is a finite end, critical and subcritical rays that
go out the end ${\sf E}$ are called {\em finite} critical or
subcritical rays; otherwise,  if the end ${\sf E}$ is an infinite
end, they are called {\em infinite}.  Furthermore, a critical or a
subcritical ray always goes out some end of $F$ .  We shall primarily be
interested in the infinite critical and subcritical rays on $F$.

Given a flute surface $F\in \mathcal{  F},$ and  an integer $n\geq 0,$ let $\alpha_n$
denote the simple closed geodesic on $F$ in the free homotopy
class of the curve $$t\rightarrow (e^{it} , n+\frac{1}{2}), 0\leq
t\leq 2\pi.$$ We shall refer to a geodesic $\alpha_n$ as a {\em dividing
loop} on  $F.$   Note that these  are well defined with the possible
exception of $\alpha_0$, which exists only if the end
corresponding to the ideal boundary $S^1\times \{0\}$ is a hole.
In what follows we shall always take this end  to be a hole
so that $\alpha_0$ does exist.

Let $\beta^* = \{(-1,t)\mid t\in (0,\infty)\}$ and  for integers $n\geq 0$ let $\gamma_n =
\{(1,t)\mid n < t< n+1\}$. Suppose in addition  that $F$ is an
untwisted flute. Then $\beta^*$ and $\gamma_n$ are geodesics since
they are fixed by the isometry $r.$ We shall refer to the geodesics $\gamma_n$ as the $\gamma$-curves of $F$.
For each $n$ the
geodesics $\beta^*$ and $\gamma_n$ are both orthogonal to the dividing loop $\alpha_n$.
Also, for any point $p$ on $\beta^*$ the geodesic ray beginning at
$p$ going out the infinite end of $F$ along  $\beta^*$ is a
critical ray. We shall refer to $\beta^*$ as the canonical
Dirichlet geodesic on $F$ and assume it to be oriented out the infinite end, with 
$\beta^*(0)\in \alpha_0.$

\vskip .1in
 \noindent {\em Proof of Proposition \ref{planar}}\, 
Suppose $S$ is a domain in the complex plane with a single infinite end.  We refer to the connected components in the complement of $S$ as complementary components. Let $\Delta_{\infty}$ denote the complementary component corresponding to the infinite end of $S$.  $S$ is endowed with the unique hyperbolic metric in its conformal equivalence class. We shall define a sequence of simple closed geodesics on $S$ so that each component in the complement of this set of geodesics on $S$ is a triply connected domain, referred to as a pair of pants, (see \cite{bas2}). 
  
Let $\Delta$ be a complementary component not containing the point at infinity.  Define the distance between $\Delta$ and $\Delta_{\infty}$ as the infimum of the (Euclidean) distances between points in $\Delta$ and $\Delta_{\infty}$ and denote this distance by $d(\Delta, \Delta_{\infty}).$ It is possible to index the complementary components not containing the point at infinity by $\NN$, so that $d(\Delta_i , \Delta_{\infty})\geq d(\Delta_{i+1} , \Delta_{\infty}).$
Let $\overline{\alpha}_1$ be a simple closed curve on $S$ that divides $\CC$ into two pieces, one of which contains only the two complementary components
$\Delta_1$ and the component  containing the point at infinity. Let $\alpha_1$ be the geodesic in the free homotopy class of $\overline{\alpha}_1$ on $S$.  Suppose the geodesics $\alpha_1,\ldots,\alpha_n$ have been defined. Let $\overline{\alpha}_{n+1}$ be a simple closed curve, disjoint from $\alpha_n$, so that the region of $\CC$ bounded by $\alpha_n$ and $\overline{\alpha}_{n+1}$ contains the single complementary component $\Delta_{n+1}$. Let $\alpha_{n+1}$ be the geodesic freely homotopic to $\overline{\alpha}_{n+1}$ on $S$. 

Let $Q$ be the set of geodesics $\alpha_i,\, i\in\NN$ defined above. Let $P_i$ denote the connected component of $S\setminus Q$ whose boundary meets the boundary of $\Delta_i$. Each of the $P_i$ is a pair of pants. Now $S$ can be reconstructed from the sets $P_i$ by \lq gluing' $P_i$ to $P_{i+1}$ along their common geodesic boundary $\alpha_i$ to get a flute surface, as in \cite{bas1}. 

To prove the converse, we simply observe that every closed curve on a flute surface $F$ divides. It is known, (see \cite{ahlfors}), that $F$ is then conformally equivalent to a domain in the plane.
  Since $F$ has a single infinite end, the proof is complete. 
 $ \hfill\Box $\\
 
%gluing%%%%%%%%%%%%%%%%%%%%%%%%%%%%%%%%%%%%%%%%%

\section{Gluing   untwisted flutes
 }\label{quilts}
Flute surfaces on which
the asymptotic geometry displays more diverse
behavior than that exhibited by untwisted flutes can be constructed
by gluing together untwisted flutes of the first kind
that have been sliced open along their canonical Dirichlet geodesics.
We shall describe a way to perform the construction to allow for infinitely many gluings
along a superstructure of scaffolding curves defined by choosing a closed subset of an oriented circle of a given circumfrence. 

\subsection{The finite steps}\label{finite}

Let $A$ denote the hyperbolic cylinder with  the  oriented
simple closed geodesic $\alpha_0$ dividing  $A$ into subsets $A^+$ and
$A^-,$ where $A^+$ is to the right of $\alpha_0.$ $A$ is completely determined by the length of $\alpha_0,$ which we denote by $a=|\alpha_0|.$ Let $C$ be a closed subset of $\alpha_0$ and let $p$ be a 
distinguished point on $\alpha_0$. The complement of $C$ in $\alpha_0$ is a countable
union of open geodesic segments which we refer to as intervals.
Order the intervals lexicographically in terms of length (larger lengths precede smaller lengths) and
oriented distance from $p,$ to get a sequence of oriented intervals
$\{I_i\}_{i=1}^{l},$ where $l,$ which may be infinity, is the number of components in the complement of $C$ on $\alpha_0.$  Henceforth, we assume that $2\leq l\leq \infty.$ Label the endpoints of $I_i,\, e_i^1 $
and $e_i^2,$ where $I_i$ is oriented from $e_i^1 $
to $e_i^2.$ Through each point $e_i^j$ there is a unique
biinfinite geodesic $\epsilon_i^j$ orthogonal to $\alpha_0$ and
oriented in the direction of the end  $A^+$, with $\epsilon_i^j(0)=e_i^j.$ Except in the case where two geodesics with different names coincide, these geodesics are
pairwise disjoint and the two ends of each $\epsilon_i^j$ go out
the two ends of the cylinder $A$. We shall refer to the the union of
the geodesics $\epsilon_i^j$ and the intervals $I_i$ as the {\it scaffolding}.
Note that the scaffolding is completely determined by the choice of the point $p,$ the orientation on $\alpha_0$ the length $|\alpha_0|$ and the set $C$. 

Let
$E_i$ denote the hyperbolic strip in $A$ bounded by the geodesics
$\epsilon_i^1$ and $\epsilon_i^2$ and containing the interval
$I_i.$
%\medskip
The scaffolding structure will serve as a foundation for the construction  in
which the hyperbolic strip $E_i$   shall be removed and then
replaced by an untwisted flute surface which has been sliced open
along its canonical Dirichlet geodesic. Consequently, in addition to
the foundational information provided, the construction also
requires a description of the flute surfaces that are to be glued
in. For each $i\in \Bbb{N},$ let $F_i$ be an untwisted flute.   We shall henceforth assume that the infinite end
on each of the untwisted flutes $F_i$ is of the 1st kind. 
One can provide sufficient conditions which guarantee this; for example, if the dividing loops grow sufficiently slowly, then the end is of the 1st kind.    

Denote the sequence of dividing geodesics on
$F_i$ by $\alpha_{i,k}$ where $k=0,1,2,\ldots$ We suppose that the length  $|\alpha_{i,0}|$ of 
$\alpha_{i,0}$
is equal to the length $|I_i|$ of the interval $I_i$. 
 On
the flute $F_i$ orient  the canonical Dirichlet geodesic 
$\beta_i^*$
in the direction of the infinite end of $F_i $ and so that $\beta_i^*(0)\in \alpha_{i,0}.$  Orient the geodesic
$\aal_{i,0}$ on $F_i$ so that the infinite end is to its right. Cut
open $F_i$  along $\beta_i^*$ to get the complete hyperbolic surface $F_i^*$ with boundary.  $F_i^*$  has the two
boundary geodesics $\beta_i^1$ and $\beta_i^2$, where  the cut open $\aal_{i,0}$ is
oriented from the the point $  \beta_i^1\cap\aal_{i,0}$  to
$  \beta_i^2\cap\aal_{i,0}$. Each $\beta_i^j$ inherits a parameterization from $\beta_i^*.$ By re-identifying the boundary
geodesics $\beta_i^1$ and $\beta_i^2$ so that $\beta_i^1(t)$ is glued to
$\beta_i^2(t),$ for $t\in \RR,$ the starting flute   $F_i$ will be reproduced.

On $F^*_i$ we shall continue to refer to the cut open dividing loops $\alpha_{i,k}$ by the 
same names. Observe that the involution $r$ of $F_i$ induces an isometric involution $r^*$ of $F_i^*$, interchanging 
$\beta_i^1(t)$ and $\beta_i^2(t)$ for $t\in \RR$, fixing the $\gamma$-curves and mapping each $\alpha_{i,k}$ onto itself.

%\medskip

With the data  consisting of $A,\, p,\, C$ and the surfaces $F_i$, we define a
sequence of flute surfaces $S_n,$ constructed by a succession of
single replacements of the type described at the beginning of the section, as follows.
Set $S_0=A.$  
Suppose that the hyperbolic surface $S_i$  has been
defined for some $i\geq0,$ and has the property that $S_{i}$ contains an
isometric copy of the original scaffolding, so that the union of
the scaffolding and the strips $E_k$  with $k>i$ embeds
isometrically in $A$. (Here objects with the same name are identified by the
embedding.)

In order to construct $S_{i+1}$ begin by removing  the interior of the strip $E_{i+1}$ from
$S_i$. Insert  $F_{i+1}^*$ in
place of  $E_{i+1}$ by identifying the two pairs of geodesics
$\epsilon_{i+1}^j$ and $\beta_{i+1}^j,\, j=1,2 ,$  so that $\epsilon_{i+1}^j(t)$ is identified with $\beta_{i+1}^j(t)$ for $t\in \RR.$
 On $S_{i+1}$ the identified pair
$\epsilon_{i+1}^j$ and $\beta_{i+1}^j $ shall be denoted by
$\beta_{i+1}^j $, which   inherits their parameterization. 
The surface $S_{i+1}$ has a hyperbolic metric defined locally at every point and is easily 
seen to be complete. 
Therefore $S_{i+1}$ is a hyperbolic surface. Furthermore, it naturally contains
a copy  of   $S_{i}\setminus E_{i+1}$, as well as a copy of  $F_{i+1}^*$.

The set of all the geodesics  $\beta_{m}^j $ on $S_n$
where $m=1,\ldots, n$ and  $ j=1,2 $ shall be denoted by $B_n$. 
%For later convenience also put $B=\bigcup_{i=1}^\infty B_i.$
 
\begin{prop}\label{prop} The surfaces  $S_n$ are all flute surfaces.
\end{prop}
\begin{proof}
First observe that if two planar surfaces are glued together in the
plane then the result is a planar surface. Each flute $F_i$, as well as its cut open relative
$F_i^*$, is planar. By induction it follows that each  of the surfaces $S_n$ is planar

Next we argue by induction that each $S_n$ has a single infinite
end. Suppose this is true for  $S_i$. Let $K^*_j$ be  a compact
exhaustion of  $S_i$ with the property that $I_{i+1} \in K^*_j$ for
each $j$ and similarly choose a  compact exhaustion
$K^{\#}_j$ of  $F_{i+1}$ so that  $\alpha_{i+1,0}\subset K^{\#}_j$ for each $j$. The $K_j^*$ correspond to the sets on $S_i\setminus E_{i+1},$ also called $K_j^*,$ which give a compact exhaustion of  $S_i\setminus\ E_{i+1}.$  Similarly, the $ K^{\#}_j$ correspond to sets on $F_{i+1}^*,$ also called $ K^{\#}_j,$ which induce a compact exhaustion of $F_{i+1}^*$. Then
$K_j=K^{\#}_j\cup K^*_j$ is a compact exhaustion of  $S_{i+1}$.
Let ${\sf E}_j$ be a nested sequence of complementary components of the sequence
$\{K_j\}.$ Then one of the following three possibilities holds: for
$j$ sufficiently large,  either  ${\sf E}_j\subset S_i \setminus E_{i+1}$ or  ${\sf E}_j\subset F_{i+1}^*$, or  ${\sf E}_j$ has
non-empty intersection with both $S_i \setminus E_{i+1}$  and $F_{i+1}^*$ for all $j$. In the first two
cases the corresponding end must be finite. In the last case it is
possible that the sets  ${\sf E}_j$ are all cylinders lying  to the left
of $\alpha_0$ and consequently,   the corresponding end is
finite. The remaining case occurs when each  ${\sf E}_j$ is a union of
open sets on $S_i \setminus E_{i+1}$ and  $F_{i+1}^*$ belonging to the infinite end of
each of those surfaces. In only this last case is it possible for
the   corresponding end to be infinite. Thus $S_{i+1}$ has a single infinite end and $S_n$ is a flute surface 
for all integers $n\geq 0$.
\end{proof}

%convergence%%%%%%%%%%%%%%%%%%%%%%%%%%%%%%%%%%%%%%%

\subsection{Geometric convergence of flute surfaces} 
Here we show that it is possible to define a flute surface as a limit of the finite constructions described in the previous section. We begin with the Fuchsian groups representing the surfaces and show that there is a way to normalize the groups so that one representing a surface $S_n$ contains the group for the surface $S_m$ when $n>m$. This leads to the definition of the surface $S_{\infty}$. Then, on the level of the intrinsic geometry, we prove that it is possible to view an arbitrarily large chunk of $S_{\infty}$ as the isometric images of a large chunk  of $S_n,$ for $n$ sufficiently large. This last fact leads to a proof that $S_{\infty}$ is a flute surface. These results will allow us to assume, in the proofs of several of the theorems that follow, that we are working on one of the surfaces $S_n,$ rather than on the more complex surface $S_{\infty}.$

\begin{lemma}\label{convergeL}
There exists a nested sequence of Fuchsian groups $\Gamma_0\subset\Gamma_1\subset\ldots, $
so that for each integer $n\geq 0$, $\HH^2/ \Gamma_n=S_n.$ Furthermore, the image of the imaginary axis $I\subset \HH^2$ covers the geodesic $\alpha_0$ on $S_n$ and the left half-plane $H^-$ covers the annular region $A^-$ on $S_n$.
\end{lemma}

\begin{proof}
Given $|\alpha_0|=a$, let $\Gamma_0$ be the group generated by the M\"obius transformation $g_0(z)=e^az.$  Observe that the imaginary axis, denoted by $I$, projects to the geodesic $\alpha_0$ on $A=\HH^2/\Gamma_0$ and the left and right half-planes project, respectively, to the subannuli $A^-$ and $A^+$.  

Define the sequence of Fuchsian groups recursively. Suppose the groups $\Gamma_0\subset\ldots\subset\Gamma_n$ have been defined. Set $E_{n+1}^+= E_{n+1}\cap A^+.$   Remove  the subsurface $E_{n+1}^+$ from $S_n$ to get the surface $S'_n.$ A loop on $S'_n$ is homotopically trivial on $S'_n$ if and only if it is homotopically trivial on $S_n.$
Therefore, the preimage $\tilde{S}'_n$ of $S'_n$ under the covering projection $\pi_n:\HH^2\rightarrow S_n$ is a connected and simply connected, $\Gamma_n$-invariant set and the restriction $\pi_n:\tilde{S}'_n\rightarrow S'_n$ is the universal covering. 

Suppose $G_{n+1}$ is a Fuchsian group representing $S_{n+1}$. Let $F_{n+1}^+= F^*_{n+1}\cap A^+.$ Then remove $F_{n+1}^+$ from $S_{n+1}$ to get the surface $S^+_{n+1}.$ Let $\tilde{S}^+_{n+1}$ be a connected preimage of  $S^+_{n+1} $ under the covering projection $\pi_{n+1}:\HH^2\rightarrow S_{n+1}.$
Then $\tilde{S}^+_{n+1}$ is simply connected and the restriction $\pi_{n+1}:\tilde{S}^+_{n+1}\rightarrow S^+_{n+1}$ is the universal covering. Let $G^+_{n+1}$ be the stabilizer of  $\tilde{S}^+_{n+1}$ in $G_{n+1}.$

By the construction, $S'_n$ and $S^+_{n+1} $ are identical. Thus there is an isometric bijection 
$\varphi:S'_n\rightarrow S^+_{n+1} $ which lifts to an isometry $\tilde{\varphi}:\tilde{S}'_n\rightarrow
\tilde{S}^+_{n+1}$. Since $\tilde{S}'_n$ and $\tilde{S}^+_{n+1}$ are subsets of hyperbolic space,
$\tilde{\varphi}$ must be the restriction of a M\"obius transformation and 
$\Gamma_n=\tilde{\varphi}^{-1}G^+_{n+1}\tilde{\varphi}.$
It follows also that $\Gamma_{n+1}=\tilde{\varphi}^{-1}G_{n+1}\tilde{\varphi} $ is a Fuchsian group representing $S_{n+1} $ and $\Gamma_{n+1}\supset\Gamma_n.$

Note that $\varphi$ maps the geodesic $\alpha_0$ on  $S'_n$ to its counterpart on $S^+_{n+1} $ and identifies the corresponding copies of $A^-$. Thus $\tilde{\varphi}$ will map $I$ to a geodesic in  $\tilde{S}^+_{n+1}$ that covers
$\alpha_0$ and takes the left half-plane $H^-$ in $\HH^2$ to a corresponding hyperbolic half-plane in $\tilde{S}^+_{n+1}$ that covers $A^-.$ It follows that for all integers $n\geq 0$, $H^-$ is precisely invariant in $\Gamma_n$ under $\Gamma_0$; that is,  $h(H^-)\cap H^- \not =\emptyset$ if and only if $h\in \Gamma_0.$ Thus, $H^-$ projects to $A^-$  and $I$ projects to $\alpha_0$ on each $S_n.$

\end{proof}

We are now in a position to define the limiting surface $S_{\infty}.$  First, define $\d{\Gamma_{\infty}=\cup_{n=0}^{l}\Gamma_n}$, where $2< l\leq\infty$   is the number of components in the complement of $C$ on $\alpha_0.$
None of the $\Gamma_n$ contain elliptic elements, and therefore the same must be true of $\Gamma_{\infty}$. Since,  for $n> 0$ $\Gamma_n$ is non-elementary, it follows that $\Gamma_{\infty}$ is   a Fuchsian group, \cite{beardon}. Define $ S_{\infty}=\HH^2/\Gamma_{\infty}.$  Let $\pi_{\infty}\rightarrow S_{\infty}$ denote the quotient map.  Set the notation $S_{\infty}=S(\alpha_0, p, C,\{F_i\})$ to emphasize the fact that the orientation on $\alpha_0,$ the point $p,$ the length $a=|\alpha_0|$, the closed set $C,$ and the sequence of untwisted flutes
$F_i$ together define the surface.  

As usual $B(q,r)\subset \HH^2$ shall denote the   set of points of distance less than r from q and we let $B^c(q,r)$ denote its closure.  
 Let $i$ be the imaginary value $\sqrt{-1}.$  Without loss of generality, we can suppose that $i$ projects to the point $p$ on $A$. Consider the translates of   the closed ball $B^c(i,r)$ by the transformations in $\Gamma_n$. Given $r>0$ there is a value  $M_r\in \NN$ so that for $ M_r\leq n\leq \infty$ and $g\in \Gamma_n,$ if $  g(B^c(i,r))\cap B^c(i,r)\not =\emptyset$ then $g\in \Gamma_{M_r}.$ This must be so since, in a Fuchsian group, only finitely many $\Gamma_{\infty}$-translates of the compact set $B^c(i,r)$ may intersect $B^c(i,r)$.

Thus,  for $n\geq M_r$, the projections $\pi_n(B^c(i,r))$ are  isometrically equivalent in the most natural   way. Note $\pi_n(B^c(i,r))$ is  connected and since  the set of transformations that create intersections is finite, the complement of $\pi_n(B^c(i,r))$ in $S_n$ has $k<\infty$ components.
Let $U_{n,1},\ldots U_{n,k-1}$ denote the complementary components of $\pi_n(B^c(i,r))$ in $S_n$ that do not contain the infinite end. Call these the finite components. For $n\geq M_r,$ define $\d{S_n(r) =\cup_{j=1}^{k-1} U_{n,j} \cup \pi_{n}(B^c(i,r))}$ 
%We may similarly define $S_{\infty}(r).$

Let $p_n\in S_n$ be the projection $\pi_n(i)$ and similarly set $\alpha_0^n=\pi_n(I).$ We have just given distinct names to the point $p$ and the geodesic $\alpha_0.$ Intrinsically, for $n\in \NN$, $S_n(r)$ is the closed ball of radius $r$ about $p_n$ in $S_n,$ with the finite complementary components   added on. It is not yet clear that this will work for $S_{\infty}$, since we have not proved that $S_{\infty}$ is  a flute.

The next proposition asserts a kind of strong geometric convergence of the surfaces $S_n$ to $S_{\infty}$.  
 
\begin{prop}\label{embed}
Given $r>e^{\frac{a}{2}}$ there is a number $N_r\in \NN$ so that for $m>n\geq N_r$ there exists an isometric bijection $\varphi_{n,m}:S_n(r)\rightarrow S_m(r) $ with $\varphi_{n,m}(\alpha_0^n)=\alpha_0^m$ and $\varphi_{n,m}(p_n)=p_m.$   Furthermore, for  $  n\geq N_r$ there exists an  isometric embedding $\varphi_n:S_n(r)\rightarrow S_{\infty}$ with $\varphi_{n}(\alpha_0^n)=\alpha_0^{\infty}=\pi_{\infty}(\alpha_0)$ and $\varphi_{n }(p_n)=p_{\infty}=\pi_{\infty}(p).$ 
 
\end{prop}
 
A collar of width $2R$ about a simple closed geodesic $\beta$ on a hyperbolic surface $S,$ written $C_S(\beta,R),$ is the set of points on $S$ of distance less than $R$ from $\beta.$ The classical Collar Lemma implies that given $R>0$ there is an number $L>0$ so that if $|\beta|<L$ then $C_S(\beta,R)$ is isometric to the collar $C_A(\alpha,R),$ where $A$ is an annulus with core geodesic of length $|\alpha|=|\beta|.$ More precisely,  $C_S(\beta,R)$ is isometric to the collar $C_A(\alpha,R)  $ if $\sinh R= (2\sinh |\beta|/2)^{-1},$ where $|\alpha|=|\beta|,$ \cite{matelski}.

One can analogously define the collars $C_{E_j }(I_j,R)$ and $C_{ F^*_j }(\alpha_{j,0}, R)$ about the geodesic segments on the  surfaces $E_j$ and $F^*_j$, respectively. Observe that identifying the boundary geodesics $\epsilon_j^1$ and $\epsilon_j^2$ of $E$ will produce an annulus. The Collar Lemma, applied to this annulus and the flute $F_j$,  implies that given $R>0$ there is a $J\in \NN$ so that if $j> J$ then $C_{E_j}(I_j,R)$ is isometrically equivalent to 
$C_{F^*_j}(\alpha_{j,0}, R)$ by an isometry taking $\alpha_{j,0}$ to $I_j.$

\begin{proof} The proposition is trivial if the number of complementary components, $l$, of $C$ in $\alpha_0$ is finite. We therefore suppose that $l=\infty.$

We normalize so that the surfaces are represented by Fuchsian groups, as in Lemma \ref{convergeL}.  Following up on the comments  preceding the statement of the proposition, given $r>0$ there is an $M_r\in \NN$ so that if $g\in \Gamma_{n}$ for $M_r\geq n\geq\infty$, and  $g(B^c(i,r))\cap B^c(i,r)\not = \emptyset,$ then $g\in \Gamma_{M_r}.$ The closed ball $B^c(i,r)$ projects to a set $B^c_n(r)\subset S_n$. Since $r>e^{\frac{a}{2}}$, $ B^c(i,r)$
contains a segment of $I$ which projects to $\alpha_0.$ 
Than for $ M_r\leq n< m\leq \infty$ there is an isometric bijection  $\varphi_{n,m}:B^c_n(r)\rightarrow B^c_m(r),$ which maps $\alpha^n_0$ to $\alpha^m_0$ and  $p_n$ to $p_m$. We suppose that the finite complementary components    $U_{n,j}$ and
$U_{m,j}$ share the corresponding boundaries on $B^c_n(r)$ and $B^c_m(r),$ respectively. We shall prove that there is a $N_r\geq M_r$ so that for any $ n\geq N_r,$ the process of creating $S_{n+1}$ from $S_n$ by excising $E_{n+1}$ and gluing in $F_{n+1}^*$, does not change the subsurface $S_n(r)$.

Suppose $n>M_r$. Since  $B^c_n(r)$ is compact, by the Collar Lemma there must exist a value $R>0$ and an integer $J_r\geq M_r$ so that for $j\geq J_r,\,  B^c_j(r)\cap  E_j\subset C_{E_j}(I_j,R).$
Consequently, the region $E_j\setminus C_{E_j}(I_j,R)$ must lie in the infinite component in the complement of $B^c_j(r)$. 

Also, by the Collar Lemma we may choose $N_r\geq J_r$ so that for $j\geq N_r$ the collar $C_{F_j^*}(\alpha_{j,0},R)$ on $F_j^*$ is isometrically equivalent to $C_{E_j}(I_j,R).$ Thus, for $n\geq N_r$, $S_{n+1}$ can be constructed from $S_n$ by replacing $ E_{n+1}\setminus C_{E_{n+1}}(I_{n+1},R) $ by $F_{n+1}^*\setminus C_{F^*_{n+1}}(\alpha_{n+1,0},R).$ Consequently, $S_{n+1}$ differs from $S_n$ only in the infinite component in the complement of $B_{n+1}^c(r)\simeq B_n^c(r).$  The first part of the proposition follows.

Let $\tilde{S}(r)$ be the connected 
$\pi_{N_r}$-preimage of $S_{N_r}(r)$, which contains $i$ and let $\Gamma\subset \Gamma_{N_r}$ denote the stabilizer of $\tilde{S}(r)$ in $\Gamma_{N_r}$. By what we have shown, it follows that $\Gamma$ is the stabilizer of $\tilde{S}(r)$ in $\Gamma_n$ for all $n\geq N_r$ and therefore in $\Gamma_{\infty}$, as well. 
The final statement of the proposition follows.
\end{proof}

Define $S_{\infty}(r)= \varphi_{n}(S_n(r))$ for some $n>N_r$. The sequence of closed balls   $\pi_{\infty}({B^c(i,k)})$, for $k\in\NN$, is a compact exhaustion of $S_{\infty}.$ Associated to this exhaustion there is an infinite end ${\sf E}$ with ${\sf E}_k=S_{\infty}\setminus  S_{\infty}(k).$ It follows from the proposition that every other end is finite. Thus we have

\begin{cor}
$S_{\infty}$ is a flute.
\end{cor}

To simplify notation, we shall henceforth dispense with superscripts and refer to the geodesic $\alpha_0$ and the point $p$ on the surfaces $S_n$  and $S_{\infty}.$

%sigma%%%%%%%%%%%%%%%%%%%%%%%%%%%%%%%%%%%%%

\section{A special class of  critical rays on $S_\infty$   }\label{number4} 

Given a point $c \in C,$ let $\sigma_c$ denote the geodesic ray on  $S_{\infty}=S(\alpha_0, p, C,\{F_i\})$ beginning at $c,$ orthogonal to $\alpha_0$ and oriented out the infinite end of
$S_{\infty}$. If $c$ lies on the boundary of one of the intervals $I_j$ in the complement of $C$, then the ray $\sigma_c$ will be contained in  the scaffolding geodesic on $S_{\infty}$ passing through $c$.

\begin{thm} \label{sigma}

For each  $c\in C,\,\sigma_c  $ is an infinite critical ray.
\end{thm}

We begin by proving several lemmas, which, along with the lemmas of the next section, shall be of use throughout the paper.  Theorem \ref{sigma} will follow as an immediate corollary of Lemma \ref{nnew5} in which we prove something slightly stronger about the ray $\sigma_c,$ namely, subarcs of $\sigma_c$ realize the distance between any point on $\sigma_c\backslash\{\sigma_c(0)\}$, and the curve $\alpha_0.$

Many of the lemmas in this and the next section will share the common setup stated below, which will be invoked repeatedly.

Let $S^*$ be either the closed set $F^*,$ an untwisted flute cut open along its canonical Dirichlet ray $\beta^*,$ or the hyperbolic strip $E^*\subset \HH^2$ bounded by two parallel, non asymptotic geodesics.  Let $\alpha_0$ denote the first, cut-open dividing loop on $F^*,$ or the common  orthogonal between the boundary geodesics on $E^*$.  The geodesics $\beta_1$ and $\beta_2$ on the boundary of $S^*$ are parameterized so that $\beta_1(0)$ and $\beta_2(0)$ are the end points of $\alpha_0$ and the geodesics $\beta_i(t),$ $i=1,2,$ go out the same end of the surface as $t \to \infty;$  in the case $S^*=F^*,$ this end is the infinite end of $F^*$.  

Let $F$ be the untwisted flute with canonical Dirichlet ray $\beta^*$ which, when cut open along $\beta^*,$ produces $F^*.$  Note that $\beta_1$ and $\beta_2$ inherit their parametrization from $\beta^*$ and that $\alpha_0$ is orthogonal to both of these boundary geodesics.  Let the union of the $\gamma$-curves of $F$ be denoted by $\Gamma.$  The curves in $\Gamma$ along with $\beta^*$ comprise the fixed-point set of the canonical isometric involution $r:F\rightarrow F.$     The involution $r$ of $F$ induces an isometric involution $r^*:F^*\rightarrow F^*$  that interchanges the boundary curves $\beta_1$ and $\beta_2.$  The fixed-point set of $r^*$  is now the set of curves (corresponding to) $\Gamma$.    Moreover, as $r^*$ is an isometry, $r^*(\beta_1(t))=\beta_2(t^\prime)$ if and only if $t=t^\prime$.

\begin{lemma}\label{nnnew1}\label{nnew3}
Let $P$ be a point in $E^*$ that does not lie on $\alpha_0.$  Then there is a unique point $Q$ on $\alpha_0$ such that the the geodesic segment $\delta$ joining $P$ to $Q$ realizes the distance $d(P,\alpha_0).$  Further, $\delta$ meets $\alpha_0$ at a right angle.  In particular, if $P$ lies on $\beta_1$ (respectively, $\beta_2$) then $Q$ lies on $\beta_1$ (respectively, $\beta_2$).
\end{lemma}

\begin{proof}
Let $\alpha$ be the full, bi-infinite geodesic containing $\alpha_0,$ and let $Q$ be the point on $\alpha$ for which the length of the geodesic arc $\delta,$ joining $P$ to $Q,$ is a minimum.  From elementary geometry, the angle at $Q$ where $\delta$ meets $\alpha$ must be a right angle.
Suppose $Q$ lies outside $E^*,$ that is, $Q$ is a point on $\alpha \backslash \alpha_0.$  Then $\delta$ must cross a boundary curve of $E^*,$ without loss of generality, say $\beta_1,$ where $\beta_1$ meets $\alpha$ at a point $R.$  Then the triangle, $\Delta PQR,$ formed by segments $PQ,$ $QR,$ and $RP$ has a right angle at $Q$ and the angle at $R$ is greater than $\pi/2.$  This is impossible.  Therefore, $Q$ lies on $\alpha_0.$ 
\end{proof}

Let $\beta$ be one of the two boundary curves of $F^*.$ 
  
  \begin{lemma}\label{notquite}\label{nnew4}
  Suppose $\alpha:[0,\tau^*]\rightarrow F^*$ is a geodesic which does not lie on $\beta$, with $\alpha(0)\in\alpha_0$ and $\alpha(\tau^*)=\beta(t^*)$ for some $t^*>0.$ Then $t^*<\tau^*$.
  \end{lemma}
  
  \begin{proof}
Consider the region $\Delta$ on $F^*$ bounded by $\alpha$, the arc $\beta([0,t^*])$  and the arc of $\alpha_0$ from $\alpha(0)$ to $\beta(0).$ If $\Delta$ is simply connected then it is a hyperbolic triangle. In that case, since $\beta$ is orthogonal to $\alpha_0$ it is a right triangle with hypotenuse $\alpha$. It follows that $t^*<\tau^*$.

Suppose now that $\Delta $ is not simply connected. The $\gamma$-curves of $F^*$ divide the surface into simply connected subsurfaces $S_1$ and $S_2$. Suppose $\beta\subset S_1.$ Use the involution $r^*$ to reflect each arc of $\alpha$ in $S_2$ to an arc on $S_1$. The arcs of $\alpha$ and $r^*(\alpha)$ on $S_1$ form a piecewise geodesic $\overline{\alpha}$, which inherits its parametrization from $\alpha$. The geodesic $\alpha'$ from $\overline{\alpha}(0)$ to $\alpha(\tau^*)=\beta(t^*)$ on $S_1$ has length less than $\alpha$. Applying the first case to  $\alpha'$ shows that its length is greater than $t^*$. Therefore, $|\alpha|=\tau^*>|\alpha'|>t^*.$  The lemma is proved.
\end{proof}
 Theorem \ref{sigma} is an immediate consequence of the following lemma.  Let $S$ be one of $S_\infty$ or $S_n$ for some $n \in \Bbb{N}.$
 \begin{lemma}\label{nnew5}
 Let $c \in C$ and $t>0.$  Then the geodesic arc $\sigma_c([0,t])$ is the unique arc that realizes the distance between $\sigma_c(t)$ and $\alpha_0.$  In particular, if $\delta : [0,\tau] \rightarrow S$ is a geodesic joining a point $Q \neq \sigma_c(0)$ on $\alpha_0$ to $\sigma_c(t),$ then $\tau>t.$ 
 
\end{lemma}
\begin{proof}
 Since $S$ is a complete, convex hyperbolic manifold, there is a geodesic arc from $\sigma_c(t)$ to some point  on $\alpha_0$ whose length is $d(\sigma_c(t),\alpha_0).$  Suppose that $Q$ is a point on $\alpha_0$ and that $\delta:[0,\tau]\rightarrow S$ is a geodesic arc joining $Q=\delta(0)$ to $\sigma_c(t).$ 

Suppose, for the moment, that $S=S_\infty.$  Since the geodesic arcs $\delta([0,\tau]),$ $\sigma_c([0,t])$ and  $\alpha_0$ are compact, there is some $r>0$ for which all of these arcs  lie in $B_{\infty}^c(r)\subset S_{\infty}^c(r).$  Thus, there is a positive integer $N_r$ such that, without loss of generality, we are working on the surface $S_n$ for some $n\geq N_r.$  If, on the other hand, $S=S_n$ for some $n \in \Bbb{N}$ then we are again working on a surface $S_n.$  In either case we may assume that we are working on a surface $S_n$ for some fixed $n\in \Bbb{N}.$ 

There are two cases to consider.  First we shall show that if $Q \neq \sigma_c(0),$ then $\tau>t.$  In the second case we shall show that if $Q= \sigma_c(0),$ but $\delta$ is distinct from $\sigma_c([0,t]),$ then $\tau>t.$  The conclusion will be that $\sigma_c([0,t])$ is the unique geodesic arc from $\sigma_c(t)$ to $\alpha_0$ realizing the distance from $\sigma_c(t)$ to $\alpha_0.$

Case I.  Recall that $B_n$ is the set consisting of the boundary geodesics of the cut-open untwisted flutes, $F_i^*,$ $i=1, \cdots, n,$ in $S_n.$  Let $\beta$ denote the full (bi-infinite) geodesic in $S_n$ that contains $\sigma_c$ so that for $t\geq 0,$ $\sigma_c(t)=\beta(t).$  Let $B^+_n = B_n \cup \{\beta\}.$ Then there is a largest point $\tau_0 \in [0,\tau]$ such that $\delta((0,\tau_0)) \cap B^+_n=\emptyset.$  Let $\beta_0$ be the bi-infinite geodesic in $B^+_n$ for which there is a $t_0>0$ such that $\delta(\tau_0)=\beta_0(t_0).$  (It is possible that $\tau_0=\tau$ and that $\beta_0=\beta.$)  The arcs $\beta_0([0,t_0])$ and $\delta([0,\tau_0])$ lie in exactly one cut-open untwisted flute or simply connected hyperbolic strip.  It then follows respectively from either Lemma \ref{nnew4} or Lemma \ref{nnew3} that $\tau_0>t_0.$ Suppose that $\tau_0 \neq \tau.$  Then the piecewise geodesic $\beta_0([0,t_0])$ followed by $\delta([\tau_0,\tau])$ is shorter than $\delta([0,\tau ]),$  i.e., $t_0+(\tau-\tau_0)=(t_0-\tau_0)+\tau<\tau.$  In this case there must be a smooth geodesic arc even shorter than this piecewise geodesic joining $\beta_0(0)$ to $\sigma_c(t).$  If $\tau_0=\tau$ then $\beta_0=\beta,$ $t=t_0$ and we have $t<\tau$ and $\beta([0,t])=\sigma_c([0,t])$ is shorter than $\delta([0,\tau]).$  In either case, there is a shorter path from $\sigma_c(t)$ to $\alpha_0$ than $\delta.$  We have shown that if $Q\neq\sigma_c(0)$ is a point on $\alpha_0,$ then a geodesic arc $\delta:[0,\tau]\rightarrow S$ joining $Q$ to $\sigma_c(t)$ does not realize the distance $d(\sigma_c(t),\alpha_0).$

Case II.  Suppose that $\delta:[0,\tau]\rightarrow S$ joins $Q=\delta(0)=\sigma_c(0)$ to $A=\sigma_c(t),$  and suppose that $\delta$ is distinct from $\sigma_c([0,t]).$  

Since $\delta$ is distinct from the arc $\sigma_c([0.t]),$  the angle at which $\delta$ meets $\alpha_0$ at $Q$ can not be a right angle.  It follows that, by moving a small distance from $Q=\sigma_c(0)$ along $\delta$ to a point $B,$ there is a point $E$ on $\alpha_0$ so that arcs $QB,$ $BE,$ and $EQ$ form a right triangle, with right angle at $E,$ and hypotenuse $QB.$  It then follows that the piecewise geodesic from $A$ to $B$ along $\delta,$ followed by $BE,$ is shorter than $\delta.$  Therefore $\delta$ does not realize the distance from $A=\sigma_c(t)$ to $\alpha_0.$  

It follows that the arc $\sigma_c([0,t])$ must be the unique geodesic arc realizing the distance from $\sigma_c(t)$ to $\alpha_0.$  In particular, as in Case I, it follows that if $Q \neq \sigma_c(0)$ is a point on $\alpha_0$ and $\delta:[0,\tau]\rightarrow S_\infty$ joins $Q$ to $\sigma_c(t), $ then $\tau>t.$
 \end{proof}
 
 The lemma above shows that $\sigma_c$ is a critical ray on all of the surfaces $S_n,$ $n \in \Bbb{N},$ and on $S_\infty.$ 

The last two lemmas of this section will be of use throughout the rest of the paper. 
The first is a generalization of Lemma \ref{nnew5}.  Again, let $S=S_\infty$ or $S=S_n$ for some $n \in \Bbb{N}.$

\begin{lemma}\label{forSigma} Let $\delta:[\tau,\tau']\rightarrow S$ be a geodesic arc for which there are not necessarily distinct points $c$ and $c'$ in $C$ such that $\delta(\tau) \in \sigma_c$ and $\delta(\tau')\in \sigma_{c'}.$ Suppose that $\delta$ is not a subarc of $\sigma_c$ or $\sigma_{c'}$ and let $t$ and $t'$ be reals so that $\delta(\tau)=\sigma_c(t)$ and $\delta(\tau')=\sigma_{c'}(t')$    Then, $\tau'-\tau > |t'-t|\ge t'-t.$
\end{lemma}

\begin{proof}
 The lemma holds if $t=t'$ so assume, without loss of generality, that $t<t'.$  Note that $|\delta |\geq d(\delta(\tau),\delta(\tau'))=d(\sigma_c(t),\sigma_{c'}(t')),$ and as $\delta$ is not a subarc of either $\sigma_c$ or $\sigma_{c'},$ the piecewise geodesic formed by the segments $\sigma_c([0,t])$ followed by $\delta([\tau,\tau']),$ is not smooth.  It follows that there is a shorter, smooth geodesic $\delta'$ joining the points $\sigma_c(0)$ and $\sigma_{c'}(t').$  Noting as observed above, that $$\tau'-\tau=|\delta|\geq d(\delta(\tau),\delta(\tau'))=d(\sigma_c(t),\sigma_{c'}(t')),$$ and employing Lemma \ref{nnew5}, we now have, 
\begin{align}
t+ (\tau'-\tau)&=d(\sigma_c(0),\sigma_c(t))+|\delta| \nonumber\\
&\geq d(\sigma_c(0),\sigma_c(t))+d(\sigma_c(t),\sigma_{c'}(t'))\nonumber\\
&>|\delta'|\nonumber\\
&\geq d(\sigma_c(0),\sigma_{c'}(t'))\nonumber\\
&>d(\sigma_{c'}(t'),\alpha_0)\nonumber\\
&=d(\sigma_{c'}(0),\sigma_{c'}(t'))\nonumber\\
&=t'.\nonumber
\end{align}
  Therefore, $\tau'-\tau>t'-t=|t'-t|.$
 
\end{proof}

Let $S^*$ be either  $F^*$ or $E^*$ as defined above.  Let $\beta$ be one of $\beta_1$ or $\beta_2.$

\begin{lemma}\label{secondlemma}
Let $\delta : [\tau_1, \tau_2]\rightarrow S^*$ be a non-trivial geodesic arc where $\delta$ is not a subarc of either $\beta_1$ or $\beta_2$ and for some some $t_1, t_2 >0,$  $\delta(\tau_1)=\beta_1(t_1)$,and $\delta(\tau_2)=\beta(t_2).$  Then $\tau_2-\tau_1>|t_2-t_1|.$

\end{lemma}

\begin{proof}
Though there are direct proofs, for efficiency, we proceed by using the results we have already obtained.  Note that if $S^*$ is a cut-open untwisted flute, we may reglue $S^*$ along $\beta_1$ and $\beta_2$ and obtain an untwisted flute on which we may apply Lemma \ref{forSigma} to obtain the result.  
%%%%%%%%%%%%%%%%%
%\note{Or, instead of below, the following.  Take your pick:  Suppose $S^*=E^*$ is a simply connected hyperbolic strip whose common perpendicular joining the boundary geodesics has length $a.$  Let $S_0=A$ be the hyperbolic cylinder with oriented, simple closed geodesic $\alpha_0$ of length $2a.$  Let $C=\{p,p'\}$ where $p$ is the distinguished point on $\alpha_0$ and $d(p,p')=a.$  Let $F_1,F_2$ be untwisted flutes with first dividing curves of length $a.$  For the flute $S(\alpha_0,p,C,\{F_1,F_2\}),$ the flute $S_1$ that results from the first step in the gluing process, has the hyperbolic strip $E^*$ as a subset.  On $S_1$, one may now apply lemma \ref{forSigma} where $\delta$ is a geodesic arc joining points on $\sigma_p$ and $\sigma_{p'}.$}
%%%%%%%%%%%%%%%%%%%%%%
  Suppose $S^*$ is a simply connected hyperbolic strip.  
 We proceed in a fashion identical to the proof of Lemma \ref{forSigma}.  We employ Lemma \ref{nnew3}.  Assume, without loss of generality that $t_2>t_1,$ and observe that the piecewise geodesic formed by $\beta([0,t_1])$ followed by $\delta([\tau_1,\tau_2])$ is not smooth.  Arguing exactly as before, we have 
\begin{align}
t_1+(\tau_2-\tau_1) &= d(\beta_1(0),\beta_1(t_1))+(\tau_2-\tau_1) \nonumber\\
&\geq d(\alpha_0,\beta_1(t_1))+d(\delta(\tau_1),\delta(\tau_2))\nonumber\\
&> d(\alpha_0,\beta_2(t_2))=d(\beta_2(0),\beta_2(t_2))\nonumber\\
&=t_2.\nonumber
\end{align}
 It follows $\tau_2-\tau_1>t_2-t_1=|t_2-t_1|.$ 
\end{proof}

 %Lemmas%%%%%%%%%%%%%%%%%%%%%%%%%%%%%%%%%%%%%%%%%%%%%%%%%%% 
 
 \section{Additional lemmas from hyperbolic geometry}\label{lemmas}

The lemmas in this section will be crucial in what follows and again involve the application of basic geometry in the hyperbolic plane. In Lemmas \ref{newFive} and \ref{newSix} we derive properties of geodesic arcs on cut-open, untwisted flute surfaces.

\begin{lemma}\label{triangle}
Let $A, B$ and $E$ be vertices of a hyperbolic
triangle where there is a right angle at $E$.  Let $a,b$ and
$e$ be the lengths of the sides opposite $A, B$ and
$E$ respectively.  Then,
$$0 <log({k^2+1 \over 2k})<e-b <\log k,$$ where $a=\log k$.
\end{lemma}

Note that for $a>0$ we have $\log \cosh(a) = \log({k^2+1 \over 2k}), $ where $a=\log k.$  Therefore, the inequality may also be written in the form, $0<\log \cosh (a)<e-b<a.$

\begin{proof}
The proof requires little more than a computation. By an
appropriate isometry we may take $E = i$, $B = ki$, $k>1$ and
$A = s+ti,$ where $s^2+t^2=1$ and $ s,t\in (0,1).$  By using the standard formula, $\rho(z,w)=\log\Big [\frac {|z-\bar{w}|+|z-w|}{|z-\bar{w}|-|z-w|}\Big ],$ for the hyperbolic distance between points $z$ and $w$ in $\Bbb{H}^2$, an elementary computation shows
that $a=\log k$ and that $f(k,t):=e-b=\log \Big [\frac {1+k^2+\sqrt{1+k^4+2k^2(1-2t^2)}}{2k(1+\sqrt{1-t^2)}}\Big ].$  Elementary multivariate calculus and some elementary algebra show that $\frac{\partial f}{\partial t} = \frac{1}{t}\Big [\frac {1}{\sqrt{1-t^2}}-\frac{1+k^2}{\sqrt{1+k^4+2k^2(1-2t^2)}}\Big ]>0$ for $0 < t < 1$.  Therefore, as a function of $t$, $f$ has no critical points and is increasing.  It follows that the minimum occurs at $t=0$   and is $\log({k^2+1 \over 2 k}).$  The maximum value occurs at $t=1$ and is $\log k.$  The former value is approached in the limit as $b \to \infty$ and the latter value is approached as $b \to 0.$
 
 \end{proof}
 
\begin{lemma}\label{fewercases}
Let $\gamma$ be a line in the hyperbolic plane, $B$ a point not
lying on $\gamma$.  Let $E$ be a point on $\gamma$ such that segment
$BE$ is perpendicular to $\gamma.$  Suppose that $A$ and $C$ are
points on $\gamma$ (not necessarily distinct from each other or
from $E$).  Let $k>1$ be such that $|BE|=\log k.$   Then $|AB|+|BC| > |AC|+2\log({k^2+1 \over 2k}).$
\end{lemma}

\begin{proof}
First note that  $|AB|>|AE|+\log({k^2+1 \over 2k}).$ If $A$ and $E$ are distinct then this is a direct consequence of Lemma
\ref{triangle}. If $A$ and
$E$ coincide, then $|AE|=0$, and $|AB|=|EB|=\log k$ so the
inequality  amounts to the observation that if $k>1$ then
$\log k > \log({k^2+1 \over 2k}).$  Similarly we have the
inequality $|BC|>|CE|+\log({k^2+1 \over 2k})$.   Next, observe that no matter what configuration the points $A,
E$ and $C$ have on $\gamma$, $|AE|+|EC|\geq |AC|.$  Thus  
$|AB|+|BC| >|AE|+ |CE|+2\log({k^2+1 \over
2k}) \geq |AC|+2 \log({k^2+1 \over 2k}).$

\end{proof}

Recall the standard setup from early in Section \ref{number4} where $F$ is an untwisted flute which, when cut along $\beta^*,$ produces $F^*$.  Let $R^*$ be the quotient, $F^*/ \langle r^*\rangle,$ and $\pi:F^* \longrightarrow R^*$ be the canonical projection. $R^*$ may be regarded as a region in the hyperbolic plane with geodesic and ideal boundary.  The geodesic boundary consists of the images of $\beta_1$, $\beta_2$ and the $\gamma$-curves of $\Gamma$ under $\pi.$ We set $\beta=\pi(\beta_1)=\pi(\beta_2)$ and, since there is little risk of confusion, we let $\Gamma$ denote the union of the $\gamma$-curves, or their projections, on $F,$ $F^*$ and $R^*.$ 

Let $\beta^\prime$ be one of  $\beta_1$ or   $\beta_2$ and let $\delta :[\tau_1,\tau_2]\rightarrow F^*$ be a non-trivial geodesic arc,  not contained in $\beta_1,$ with  $\delta (\tau_1)\in\beta_1$ and   $\delta (\tau_2)\in\beta^\prime.$  
Since $[\tau_1,\tau_2]$ is a closed set, there must be a point $\lambda_0 \in [\tau_1,\tau_2]$ such that $d(\delta(\lambda_0),(\beta_1 \cup \beta_2)) \ge d(\delta(\lambda),(\beta_1 \cup \beta_2))$ for all $\lambda \in [\tau_1,\tau_2].$

\begin{lemma}\label{newFive}\label{lambdaInvolution}
Let $\delta : [\tau_1,\tau_2]\longrightarrow F$, $\beta_1$ and $\beta_2$ be as above.  Then there is a point $\lambda_0 \in [\tau_1,\tau_2]$ such that $d(\delta(\lambda_0),(\beta_1 \cup \beta_2)) \ge d(\delta(\lambda),(\beta_1 \cup \beta_2))$ for all $\lambda \in [\tau_1,\tau_2],$ and further $\delta(\lambda_0)$ lies on $\Gamma$, the fixed axis of $r^*.$
\end{lemma}
\begin{proof}
It only remains to show that $\delta(\lambda_0)$ lies on  $\Gamma.$  Let $B^\prime=\delta(\lambda)$ for some $\lambda \in (\tau_1,\tau_2).$  Suppose that $B^\prime$ does not lie on a $\gamma$-curve.  Let $B=\pi(B^\prime) $ and let $\delta^*=\pi(\delta),$ which is a piecewise geodesic arc on $R^*$.  

There are two cases to consider.  
In the first case $B$ lies on a geodesic segment of $\delta^*$ joining a point $A$ on $\beta$ to a point $D$ in $  \Gamma.$ The point $D$ lies on a geodesic $\gamma\subset\Gamma$ which  is disjoint from $\beta.$    It follows from elementary hyperbolic geometry that $B$ is closer to $\beta$ than $D$; to wit:  

Let $C$ be the point on $\beta$ such that segment $DC$ is perpendicular to $\beta.$  If $C=A$ then, clearly, $|AB|<|AD|$ since $B$ lies between $A$ and $D$.  If $C \ne A$ then consider the hyperbolic triangle $\triangle ADC.$  The point $B$    lies on side $AD$ and must be closer to the line through segment $AC$ than the opposite vertex $D$.

Now we turn to the second case, in which $B$ lies on a geodesic segment $\alpha$ of $\delta^*$ joining two points on $\Gamma$.  The full geodesic containing the arc $\alpha$  meets two of the $\gamma$-curves on the boundary of $R^*$ and therefore $\alpha$ must be disjoint from $\beta$.  Let $A$ and $C$ be points on $\alpha$ so that $B$ lies between $A$ and $C$.  It is an elementary fact in hyperbolic geometry that one of $A$ or $C$ is farther from $\beta$ than $B$.

It follows that $\delta(\lambda_0)$ must lie on $\Gamma$, the fixed axis of $r^*$.

\end{proof}

Consider again the situation set up in the paragraph proceeding Lemma \ref{newFive}.  
\begin{lemma}\label{newSix}
Let $\delta :[\tau_1,\tau_2] \rightarrow F^*$ be a geodesic arc joining $\beta_1$ and $\beta'$ as above.  Let reals $t_1$ and $t_2$ be chosen so that $\delta(\tau_1)=\beta_1(t_1)$  and $\delta(\tau_2)=\beta'(t_2).$  The arc $\delta$ crosses $\Gamma $ at a point  $\delta(\lambda_0)$, for some $\lambda_0 \in [\tau_1,\tau_2].$  Let $k>1$ be chosen so that $\log k = d(\delta(\lambda_0),(\beta_1 \cup \beta_2)).$ Then $\tau_2-\tau_1 > |t_2-t_1|+2 \log\big [\frac {k^2+1}{2k}\big ].$
\end{lemma}

\begin{proof}
Let $R^*$, as above, be the planar quotient surface $F^*/\langle r^*\rangle.$  We shall employ Lemma \ref{fewercases}, where $\gamma$ in the lemma, corresponds to the boundary geodesic $\beta=\pi(\beta_1)=\pi(\beta')$,  $B=\pi(\delta(\lambda_0))$ is the projection into $R^*$ of the point at which $\delta$ crosses $\Gamma$, $A=\pi(\beta_1(t_1))=\beta(t_1)$, and $C=\pi(\beta'(t_2))=\beta(t_2)$.  

The length of $\delta$ is equal to the length of the piecewise geodesic $\delta^*:=\pi(\delta).$ Suppose $\delta $ crosses  $\Gamma$   exactly once.   Then  
$$\tau_2-\tau_1=|\delta|=|\delta^*|  = |AB|+|BC|.$$   
On the other hand, if $\delta$ crosses $\Gamma$ more than once then  
$$\tau_2-\tau_1=|\delta|=|\delta^*| \geq |AB|+|BC|.$$ 
  In either case, $\tau_2-\tau_1 \ge |AB|+|BC|.$  
  
  Let $E$ be the point on $\beta$ such that segment $BE$ is perpendicular to $\beta$ and note that $|BE|=d(\delta^*(\lambda_0),\beta)=\log k.$  The hypotheses for Lemma \ref{fewercases} are satisfied and therefore, $\tau_2-\tau_1 \ge |AB|+|BC|>|AC|+ 2 \log\big [\frac {k^2+1}{2k}\big ]=|t_2-t_1|+2 \log\big [\frac {k^2+1}{2k}\big ].$

\end{proof}

We will  need yet another  result from hyperbolic geometry.
\begin{lemma}\label{quad}
 Let $l_1$ and $l_2$ be asymptotic geodesics in $\Bbb{H}^2$ which both limit at a point $Q\in\partial\HH^2.$ Suppose  that $B\in\l_1$ and $A\in l_2$ are points so that the arc $AB$ is orthogonal to $l_1.$ Further suppose that $D$ is a point on $l_1$ between $B$ and $Q$ and $C$ is a point on $l_2$ between $A$ and $Q$ so that the arc $CD$ is also orthogonal to $l_1.$
   Then  $|AC|+|CD|<|AB|+|BD|.$

\end{lemma}

Observe that this lemma is about a special sort of quadrilateral, $BACD$, where the sides $BD$ and $AC$ lie on parallel, asymptotic geodesics. 
\begin{proof}
Let $m$ denote the full geodesic in $\HH^2$ containing the arc $AB$. Without loss of generality we may take $l_1$ to be the imaginary axis; $m$ the semicircle passing through $-1,$ $i$ and $+1$;  $B=i$; $A$ is the point $x_0+iy_0$ for some $x_0$, $0<x_0<1$; and $l_2$ is the line $Re z = x_0.$  Let $C$ be the point $x_0+iy_1$, $y_1>y_0$.  Then $D$ is the point $ri$ where $x_0^2+y_1^2=r^2$.  Using the standard formula, $\rho (z,w)=\log \big ( \frac {|z-\bar{w}|+|z-w|}{|z-\bar{w}|-|z-w|}\big )$, for the distance between two points, $z$, $w$, in the hyperbolic plane, an entirely elementary calculation shows that: 
$|BD|=\log r;$ 
$|AB|=\log \big (\frac {1+x_0} {y_0} \big );$ 
$|AC|=\log \big (\frac {y_1} {y_0}\big );$ and
$|CD|= \log \big (\frac {r+x_0} {y_1}\big).$  Therefore, 
$$|AB|+|BD|-(|AC|+|CD|)=\log (1+x_0)+\log\big (\frac {r} {r+x_0}\big ).$$
For fixed $x_0,$ the right hand side of the equation is a strictly increasing function of $r$.  Moreover, when $r=1$ the right hand side above has value $0$.  Therefore,  $|AC|+|CD|<|AB|+|BD|$.

\end{proof}

%dirichlet%%%%%%%%%%%%%%%%%%%%%%%%%%%%%%%%%%% 
  
 \section{Specifying the boundary of a Dirichlet polygon}\label{example2}

\subsection{Definitions and a theorem}\label{def}
Let $G$ be a Fuchsian group acting on the hyperbolic plane $\HH^2$ and representing the surface $S=\HH^2/G. $ Given $\tilde{p}\in\HH^2$, one defines the Dirichlet polygon $D(\tilde{p},G)=D$ of $G$ centered at $\tilde{p}$ to be the set of all points $q\in \HH^2$ so that $d(\tilde{p},q)\leq d(g(\tilde{p}),q)$ for all $g\in G\setminus\{id\}.$

 Let $p$ be  the projection of $\tilde{p}$ to $S$ and let $\sigma$ be a critical ray on $S$ with $\sigma(0)=p$. Then $\sigma$ lifts to a geodesic ray $\tilde{\sigma}$ beginning at $\tilde{p}$, which is entirely contained in $D$. Recall that the closure of a set $X$ is written $X^c$ and $\hat{\CC}$ denotes the Riemann sphere. Let $\partial_{\infty}D$ be that part of the boundary of $D^c\subset \hat{\CC}$ lying at infinity; that is, in the extended real line $\hat{\RR}$. Then we have $\d{\lim_{t\rightarrow\infty}\tilde{\sigma }(t)\in \partial_{\infty}D}.$ Conversely, if $q$ is a point in $\partial_{\infty}D$, then the ray $\tilde{\sigma}$   beginning at $\tilde{p}$ and limiting at $q$  projects to a critical ray on $S$.
These observations follow easily from the definitions of a critical ray and a Dirichlet polygon.
 Clearly, there is an intimate relationship between critical rays and the boundary points of Dirichlet polygons.

We show how quilted surfaces can be used to prove the following

\begin{thm}\label{partial}
Given any compact set $C^*\subset \RR,$ which is not an interval, there exists a Fuchsian group $\Gamma_{\infty}$, representing a quilted surface $S_{\infty}$, and a point $\tilde{p}\in\HH^2$ so that $\partial_{\infty}D(\tilde{p}, \Gamma_{\infty})$ consists of the union of $C^*$, a countable set of isolated parabolic fixed points of  $\Gamma_{\infty}$ and an interval.
\end{thm}

The interval in the 
theorem comes from the finite end of the surface $S_{\infty}$ associated with the annulus $A^-$. If we were to cut out $A^-$ and glue a twice punctured disc into its place, the boundary of the Dirichlet polygon would not contain the interval. We shall do exactly this at the end of the section, to prove Theorem \ref{C}.

The idea of the proof of Theorem \ref{partial}
is to begin with an annulus $A$ with a dividing geodesic $\alpha_0$ of length $a\leq 1$ and to choose a closed set $C$ on $\alpha_0$  that corresponds to the set $C^*$. We then choose appropriate untwisted flute surfaces $F_i$ to define the surface $S_{\infty}=S(\alpha_0, p , C, \{F_i\})$ so that for a fixed point $ p\in \alpha_0$ and for all, but possibly one value of, $ c \in C$ there exists a unique critical ray beginning at $ p $ which is asymptotic to $\sigma_c$ and, furthermore, these are all of the infinite critical rays on $S_{\infty}.$ This surface will be uniformized by a Fuchsian group for which the theorem holds

\subsection{Working on quilted surfaces}\label{R}
We begin by constructing quilted surfaces like those described above.
Suppose $|\alpha_0|=a$ and $C\subset \alpha_0$ is a closed set. Let  $\log k=R$ and define $\kappa(R)=2
\log({k^2+1 \over 2k})=2\log \cosh (R)$. Let $k>1$ be a value for which   $\kappa(R)=a+1.$
 For each $i\in \NN,$  let $F_i$ be an untwisted flute on which all of the finite ends, except the end associated to the geodesic $\alpha_{i,0}$, are punctures. Further suppose that for $j>0$ the dividing curves $\alpha_{i,j}$ have length $R$. It follows  easily from this choice, that each of the $F_i$ has an infinite end of the first kind. We shall maintain this stipulation on the structure of the $F_i$ for the remainder of the section. Also, $S_{\infty}=S(\alpha_0, p, C,\{F_i\})$ shall be a surface constructed with these assumptions. Observe that with our choice of flutes $F_i$, the surface $S_{\infty}$ is  the same, independent  of the choice of the point $p$.
 
Fix $ p \in \alpha_0.$ Given $c\in C$ define the piecewise geodesic rays $ \overline{\delta}_c^s$ and  $ \overline{\delta}_c^l$. The ray $ \overline{\delta}_c^s$ is formed by the shorter arc of $\alpha_0$ from $ p $ to $c$, followed by $\sigma_c$. $ \overline{\delta}_c^l$ is formed in a similar fashion using the longer arc of $\alpha_0$ from   $ p $ to $c$. In the special case where the point $c$ is half-way around $\alpha_0$ from $p$,  the two arcs of $\alpha_0$ are of equal length and we make an arbitrary choice of which is  $ \overline{\delta}_c^s$ and which is $ \overline{\delta}_c^l$.
Let $\delta_c^s$ and $\delta_c^l$ be the geodesic rays beginning at $ p $ which are asymptotically homotopic to $ \overline{\delta}_c^s$ and $ \overline{\delta}_c^l$, respectively. By this we mean that, for example, there are lifts of  $ \overline{\delta}_c^s$ and $\delta_c^s$ to $\HH^2$ which have the same initial point and the same endpoint at infinity. Observe that both $\delta_c^s$ and $\delta_c^l$   are asymptotic to the critical ray $\sigma_c.$

\begin{prop}\label{onepointcritical}
For each $c\in C$ one of $\delta_c^s$ or $\delta_c^l$ is a critical ray. If $\delta$ is an infinite critical ray on $S_{\infty}$ beginning at $p$, then for some $c\in C$, $\delta=\delta_c^s$ or $\delta=\delta_c^l$.
Moreover, there can be at most one value of $c$ for which both rays are critical rays.
% In the case where  $c$ is half-way around $\alpha_0$ from $p$, both  $ \delta_c^s$ and $  \delta_c^l$ are critical.
\end{prop}

As a consequence of the proof of Theorem \ref{partial}, we will be able to conclude that $\delta_c^s$ is always critical.

We begin with some notation and prove a lemma. Let $F^*$ be a cut-open, untwisted flute bounded by geodesics $\beta$ and $\beta'$. Suppose $\epsilon $ and $\delta$ are geodesics with their initial points on $\beta$ and their terminal points on $\beta'$.
We shall say that $\epsilon$ is $\beta$-{\it homotopic} to  $\delta$ if
  $\epsilon $ is  homotopic to $\delta$ by a homotopy that does not move the initial and terminal points of $\epsilon$ off the geodesics  $\beta$ and $\beta'$. 
Observe that if $\epsilon$ is $\beta$-homotopic to a dividing curve $\alpha_j$ on $F^*$ then it crosses a single $\gamma$-curve on $F^*$; that is, it crosses the unique $\gamma$-curve which is orthogonal to  $\alpha_j$.
  
  On the same surface $F^*$, recall that $\beta $ and $\beta'$ are parameterized so that $\beta(0), \beta'(0)\in\alpha_0$ and both $\beta $ and $\beta'$ go  out the infinite end as $t\rightarrow\infty.$ Also, the canonical involution $r^*$ fixes the $\gamma$-curves, leaves invariant the dividing curves and interchanges $\beta$ and $\beta'$.
 Let  $ \overline{\delta}_c$ denote one of the geodesics $ \overline{\delta}_c^s$ or $ \overline{\delta}_c^l.$

\begin{lemma}\label{plus1}
Let $\alpha:[0,\tau^*]\rightarrow S_{\infty}$ be a geodesic with $\alpha(0)=p$ and $\alpha(\tau^*)= \overline{\delta}_c(t^*)$ for some $t^*\in [0,\infty).$ Suppose  there is an arc of $\alpha$ that is $\beta$-homotopic to a dividing curve $\alpha_i\not = \alpha_0$ on the cut-open, untwisted flute $F^* \subset S_{\infty}.$ Then $\tau^*>t^*+1.$
\end{lemma}
%%%%%%%
\begin{proof}
By hypothesis there is a cut-open flute $F^*$ bounded by geodesics $\beta$ and $\beta'$ and a connected segment of $\alpha \cap F^*$ that is $\beta$-homotopic to $\alpha_i\neq \alpha_0.$  Let $0\leq \tau <\tau',$ and $t,\,t' \geq 0$ and $\beta,\,\beta'$ be such that $\alpha([\tau,\tau'])$ is $\beta$-homotopic to $\alpha_i$, $\alpha(\tau)=\beta(t),$ and $\alpha(\tau')=\beta'(t').$  Since $|\alpha_i|=R,$ by Lemma \ref{newSix}, $\tau'-\tau>|t'-t|+\kappa(R)\geq t'-t+a+1.$  It also follows from Lemma \ref{nnew5} that $\tau>t.$ 

Let $t_c$ be the length of the subarc of  $\overline{\delta}_c$ that runs along $\alpha_o$ and joins $p$ to $c.$  The length of the scaffolding curve lying along  $ \overline{\delta}_c$ from the point $c$ on $\alpha_0$ to $\alpha(\tau^*)$ is then $t^*-t_c.$ Note also that $a-t_c>0.$  Suppose $\tau^*\neq\tau'.$ Then by Lemma \ref{forSigma}, $\tau^*-\tau'>|t^*-t_c-t'|> t^*-t_c-t'.$  Putting this all together we have,
\begin{align}
\tau^*&=(\tau^*-\tau')+(\tau'-\tau)+\tau \nonumber\\
&>(t^*-t_c-t')+(t'-t+a+1)+t \nonumber\\
&=t^*+(a-t_c)+1>t^*+1. \nonumber
\end{align}

If $\tau^*=\tau'$ we also have $t^*-t_c=t'.$
Then, similar to the above, we have $\tau^*=(\tau^*-\tau)+\tau>(t^*-t_c-t+a+1)+t=t^*+(a-t_c)+1>t^*+1.$  Therefore, in either case, $\tau^*>t^*+1.$
\end{proof}

Let $\delta_c$ denote one of the geodesics $\delta_c^s$ or $\delta_c^l$ and  denote the other one by $\delta'_c$

\vskip .1in
 \noindent {\em Proof of Proposition \ref{onepointcritical}}\, 
To  begin we prove the first statement of the proposition. Let $\alpha:[0,u^*]\rightarrow S_{\infty}$ be a geodesic with  $\alpha(0)=p$ and $\alpha(u^*)=\delta_c(s^*)$ for some $u^*, s^*\in [0,\infty).$  We will prove that one of $\delta_c$ or $\delta'_c$ is critical by  showing that for any such geodesic $\alpha$, either  $s^*\leq u^*$ or $\delta'_c$ is critical. We argue by contradiction.
Suppose there exists a geodesic $\alpha$, as above, with $u^*<s^*.$ We may further assume that $\alpha$ realizes the distance between $p$ and $\alpha(u^*)$ and is therefore a minimal length geodesic between its endpoints. In particular, $\alpha$ is simple. As described earlier, there is no loss of generality in assuming that $S_{\infty}$ is one of the finitely glued surfaces $S_n.$

The first case to consider is where  $\alpha$ contains a subarc which is $\beta$-homotopic to a dividing curve $\alpha_{i.j}$ with $j\not =0$ on a cut-open subflute $F_i^*$.  
Choose a minimal length geodesic $\mu$ from $\delta_c(s^*)$ to $\overline{\delta}_c,\,$ where $\mu(0)=\delta_c(s^*)$ and $\mu(\tau^*)=\overline{\delta}_c(t_*)$ for values $t^*, \tau^*\in [0,\infty)$.
Define the piecewise geodesic arcs $\alpha^*$ and $\delta_c^*$ by adjoining the arc $\mu$ to   $\alpha$ and to the arc of $\delta_c,$ between $p$ and $\delta_c(s^*)$, respectively. In particular, 
\begin{equation} \label{alpha*}
\alpha^*(u)=\left\{
\begin{array}{ll}
\alpha(u) & 0\leq u \leq u^*\\
\mu(u-u^*) & u^* \leq u \leq \tau^*+u^*
\end{array}  \right.
\end{equation} 
and
\[
\delta^*_c(s)=\left\{
\begin{array}{ll}
\delta_c(s) & 0\leq s \leq s^*\\
\mu(s-s^*) & s^* \leq s \leq \tau^*+s^*.
\end{array}\right.
\]
Together, $\delta^*$ and $\overline{\delta}_c $ bound a hyperbolic quadrilateral  that satisfies the hypothesis of Lemma \ref{quad}. It follows that $\tau^*+s^*<t^*$ 

Let $\rho$ be the geodesic homotopic to $\alpha^*$ relative to endpoints. $\rho$ crosses the same $\gamma$-curves as $\alpha$ on subflutes. Therefore, $\rho$ contains a subarc which is $\beta$-homotopic to $\alpha_{i,j}$ on $F_i^*.$ Then employing the above inequality and Lemma \ref{plus1}, we have
$$\tau^*+u^*= |\alpha^*| > |\rho| > t^*+1> \tau^*+s^*+1.
$$
Thus, $u^*>s^*$, which gives a contradiction.  

Now  consider the case where no subarc of $\alpha$ is $\beta$-homotopic to an arc
$\alpha_{i,j}$ with $j\not =0$ on a cut-open subflute $F_i^*$.   Then, in order, $\alpha$ crosses the curves $\beta_1, \beta_2,\ldots, \beta_m\in B^+_n$ Recall that $B_n^+$ is $B_n\cup\{\beta\}$ where $\beta$ is the full geodesic containing $\sigma_c.$

First we suppose that two   consecutive curves $\beta_j$ and $\beta_{j+1}$ are equal. Then there is a  cut-open, untwisted subflute $F^*, $ one of whose boundary geodesics is $\beta_j$, so that    $\alpha\cap F^*$ contains an arc $\epsilon$ both of whose endpoints lie on  $\gamma$-curves of $F^*$. 

Apply the canonical involution $r^*$ to get the geodesic arc $r^*(\epsilon)$. Since the $\gamma$-curves of $F^*$ are fixed by $r^*$, $\epsilon$ and $r^*(\epsilon)$ share the same endpoints. Thus we may replace the   geodesic segment $\epsilon $ of $\alpha$ by the arc $r^*(\epsilon)$. This results in a new piecewise geodesic joining $p$ to $\alpha(u^*).$ Taking the geodesic freely homotopic to this curve relative to endpoints gives a geodesic $\alpha'$ which is shorter than $\alpha.$ This contradicts the assumption that the length of $\alpha$ is the distance between its endpoints. We may then suppose that all of the curve $\beta_j$ are distinct. Possibly abusing notation, let $F_j^*$ denote the cut-open flute surface bounded by $\beta_j$ and $\beta_{j+1}.$

We would like to show that $\rho =\alpha\cap F_j^*$ is $\beta$-homotopic to $\alpha_{j,0} $ for each cut-open subflute $F_j^*, j=1,\ldots m-1.$  If not, then by earlier considerations, it    cannot  be $\beta$-homotopic to $\alpha_{j,k}$ for any $k\in\NN.$ Consequently, $\rho$ must contain at least 3 intersections with $\gamma$-curves on some $F_m^*$ and there will be an arc $\epsilon$ of $\rho$ whose endpoints lie on distinct $\gamma$-curves of $F_m^*.$ As above, replace the arc $\epsilon$ of $\alpha$ by its reflection $r^*(\epsilon)$. Taking the geodesic homotopic to this piecewise geodesic path results in a shorter geodesic between the endpoints of $\alpha,$ which is a contradiction.
We may now suppose that for each of the cut open flutes   $F_j^*$, $\alpha\cap F_j^*$ is $\beta$-homotopic to $\alpha_{j,0}.$

If none of the $\beta_{j}$ contains the ray $\sigma_c$ then $\alpha$ will cross exactly the same curves in $B_n$ as the geodesic arc of $\delta_c$ from $p$ to $\delta_c(s^*)=\alpha(u^*). $  For this to happen $\alpha$ must be homotopic to the arc of $\delta_c$ and consequently it will actually coincide with that  arc.  Therefore $s^*=u^*.$
 
 If one of the geodesics $\beta_{j}$ does contain the ray $\sigma_c$ then it must be the last one,
 $\beta_{m}$.   Now consider the ray $\delta'_c.$ Either it is critical or there is a  geodesic $\alpha':[0,u']\rightarrow S_{\infty}$  with  $\alpha'(0)=p$, $\alpha(u')=\delta'_c(s')$ for some $u', s'\in [0,\infty)$  and $u'<s'.$ 
 
 Adjust the choice of $n>N_r$ so that $ B^c_n(r)\supset \alpha'$ as well. The first part of the proof will be completed by showing that both $\alpha$ and $\alpha'$ cannot realize the distance between their endpoints. As a consequence of the preceding arguments, if $\alpha'$ realizes the distance between its endpoints then, as with $\alpha$,  $\alpha'$ will cross each geodesic in $B^+_n$ at most once and the last geodesic in $B^+_n$ that $\alpha'$ crosses must contain the ray $\sigma_c$. Since $\delta_c$ and $\delta'_c$ limit at $\sigma_c$ from opposite sides, the geodesics $\alpha$ and $\alpha'$ must intersect; that is, $\alpha(d)=\alpha'(e)$ for some $d<u^*$ and $e<u'.$ Without loss of generality suppose that $d\leq e$. Then consider the piecewise geodesic arc which is $\alpha$ from $p$ to $\alpha(d)$
 followed by $\alpha'$ from $\alpha(d)=\alpha'(e)$ to $\alpha'(u')$. The geodesic freely  homotopic to this arc relative to endpoints, goes from $p$ to $\alpha'(u')$ but is shorter than $\alpha'$, showing that $\alpha'$ does not realize the distance between its endpoints. That completes the proof that at least  one of  $\delta_c^s$  or $\delta_c^l$  is critical.
 
The second point to the proposition is that the only possible infinite critical rays beginning at $p$ are the rays
$\delta_c^s$ or $\delta_c^l$ for $c\in C$. Suppose, to the contrary, that there is an infinite critical ray $\alpha$ distinct from the above ones. 
All critical rays beginning at $p$ intersect only at $p$.   Consequently, there would be   values $c_1$ and $c_2$ in $C$ bounding an interval on $\alpha_0$ so the $\alpha $ would lie in the region on $S_{\infty}$ bounded by $\delta_{c_1}$ and $\delta_{c_2}$.

Then there is a  cut-open flute subsurface $F^*$ bounded by $\beta_1\supset\sigma_{c_1}$ and $\beta_2\supset\sigma_{c_2}$. Let $\alpha_i,\, i\geq 0$ denote the dividing curves on $F^*$ and let  $\gamma_i,\, i\geq 0$ be the $\gamma$-curves on $F^*$  where $\gamma_i\cap\alpha_j \not =\emptyset$ if and only if $i=j.$ Throughout the paper we have assumed that on all flutes $F_i$, the infinite end is of the first kind. It follows from \cite{haas}, that a geodesic ray going out the infinite end of $F$ must either cross the $\gamma$-curves or be asymptotic to the canonical Dirichlet ray $\beta^*.$ Therefore, in order for $\alpha $ to be distinct from $\delta_{c_1}$ and $\delta_{c_2}$, it must intersect at least one of the $\gamma$-curves other than 
$\gamma_0$. If $\alpha$ were to intersect two of the $\gamma$-curves then, as before,  we could use the involution to produce a shorter geodesic between $p$ and a point on $\alpha$. Since $\alpha$ is critical, this is not possible. Thus $\alpha$ would have to meet exactly one curve $\gamma_k$ for $k>0$ and then, beyond that, would   be asymptotic to either $\sigma_{c_1}$ or $\sigma_{c_2}$. Without loss of generality we take it to be the former and write $c_1=c$.
Then $\alpha$ would also be eventually asymptotic to $\delta_{c}$.

Since $\alpha$ is   asymptotic to $\delta_{c},$ which is itself asymptotic to $\overline{\delta}_c$, we can find values $t^*, s^*, u^*\in [0,\infty)$ so that  the following inequalities are satisfied by the distances to 
 $\overline{\delta}_c$: $d(\alpha(u^*), \overline{\delta}_c(t^*))<1/3$ and   $d(\delta_c(s^*), \overline{\delta}_c(t^*))<1/3.$ Let $\mu $ be the geodesic from $\alpha(u^*)=\mu(0)$ to $\overline{\delta}_c(t^*)=\mu(\tau^*)$ and let $\nu$ be the geodesic from $\delta_c(s^*) =\nu(0)$ to $\overline{\delta}_c(t^*)=\nu(\zeta^*).$
 
 Define $\alpha^*$ as in Equation (\ref{alpha*}). Also, define the piecewise geodesic
 $$
 \delta_c^*(s)=\left \{
 \begin{array}{ll}
 \delta_c(s) & 0\leq s \leq s^*\\
 \nu(s-s^*) & s^* \leq s \leq s^*+\zeta^*\\
 \mu^{-1}(s-s^*-\zeta^*) & s^*+\zeta^* \leq s \leq s^*+ \zeta^*+ \tau^*,
 \end{array}
 \right.
 $$
 Where $\mu^{-1}$ is the  geodesic $\mu$ traversed in the opposite direction.
 
 Let $\rho$ and $\delta^*$ denote the geodesics homotopic to $\alpha^*$ and $\delta_c^*,$ respectively 
 Note that  $\delta^*$ is a geodesic from $\alpha(0)$ to $\alpha(u^*).$
As a consequence of  Lemma \ref{plus1}, applied to $\rho$, $u^*+\tau^* > |\rho |> t^* +1.$  By the definition of $\delta_c^*$, and considering its restriction to the interval $[0, s^*+\zeta^*],$ we have $s^*< t^* + 1/3.$ Then  
 $$
  u^*+\tau^*> |\rho |\ > t^* +1  >s^* +\zeta^*+
\tau^*>|\delta^*|
$$
which shows that $\alpha$ cannot be critical. That completes the proof of the second statement of the proposition.

Now we need to see that there is at most one $c\in C$ for which $\delta_c^s$  and  $\delta_c^l$ are both critical. If not, and there were a second $c'\in C$ so that $\delta_c'^s$  and  $\delta_c'^l$  are critical, then two of the four critical rays would have to intersect at a point other than $p$. But since critical ray cannot intersect,  this is impossible
 $ \hfill\Box $\\
\\

 Observe that by considering the ideal triangle with sides $\alpha_0, \delta_c$ and $\delta_c'$, one can show that both $\delta_c$ and $\delta_c'$ are critical if the distance from $c$ to $p$ along $\alpha_0$ is $a/2.$

We shall also need the following lemma in our proof of Theorem \ref{partial}.

\begin{lemma}\label{ends}
Let ${\sf E}$ be a finite end of the first kind on a surface $S$ and let $p$ be a point on $S$. Then there exist only finitely many critical rays beginning at $p$ going out the end ${\sf E}$.
\end{lemma}

\begin{proof}
Choose a number $M>0$ so that in the complement of the ball $B(p,M)$ there is a component $V$ which is   a punctured disc containing the end ${\sf E}$.
Let $m$ denote the length of the boundary of $B(p,M)$.

Suppose $\alpha$ is a geodesic ray going out the end ${\sf E}$ that intersects a component $U\not \subset V$ in the complement of $B(p,M+m)$. Then there is a subarc of $\alpha$ which intersects $U$ and has its endpoints in $\partial B(p,M)$. This arc of $\alpha$ has length greater than $m$. Therefore, replacing this arc by a curve in $\partial B(p,M)$ joining its endpoints, results in a shorter curve from $p$ to any point of $\alpha$ lying beyond the arc. This shows that $\alpha$ cannot be a critical ray.

Thus each critical ray from $p$ out ${\sf E}$ must lie in $V\cup B(p,M+m)$. The  end ${\sf E}$ is of the first kind. Therefore, there cannot exist two geodesic rays beginning at $p$, that go out {\sf E} and bound a simply connected region on $S$. Since $B(p,M+m)$ has finitely many complementary components, any set of simple disjoint geodesic rays beginning at $p$ and going out $E$ must be finite. In particular, there can only be finitely many critical rays.
\end{proof}

\subsection{The Dirichlet polygon}

Let $A$ be the hyperbolic cylinder with oriented dividing geodesic $\alpha_0$ of length $a\leq1.$  $A$ can be uniformized by a Fuchsian group $\Gamma_0,$  generated by the transformation $g(z)=e^az.$ The covering projection $\pi:\HH\rightarrow A$ takes the imaginary axis $I$, oriented from 0 to $\infty,$ to the oriented geodesic $\alpha_0.$ With this setup, the left half-plane $H^-$ covers the subannulus $A^-$ and the right half-plane covers the subannulus $A^+.$ We fix the point $\tilde{p}=ie^{\frac{a}{2}}\in I$ and let $p=\pi(\tilde{p})\in \alpha_0.$

Given a compact set $C^*\subset\RR$, we shall define a projection of $C^*$ to a closed set $C$ on $\alpha_0$, which   shall be used to construct a quilted surface. There is a M\"obius transformation $\varphi(z)=Az+B,\, A, B \in \RR, A>0, $  taking $C^*$ into the interval $[1,e^a]$ so that $1,e^a \in\varphi( C^*)$.  

Given $x\not = 0,$ let $\psi(x)$ denote the hyperbolic geodesic in $\HH^2$ with endpoints $x$ and $-x$. Define the map $\Lambda:C^*\rightarrow I$, that takes $c^*\in C^*$ to the point $I\cap \psi(\varphi(c^*)).$ Define $\tilde{C}=\Lambda(C^*).$ Then
$\pi\circ\Lambda:C^*\rightarrow \alpha_0$ defines a map which is one-to-one, except for identifying the endpoints of $C^*.$ Define $C=\pi\circ\Lambda(C^*)$. 

As in Section \ref{finite}, $C$ defines a sequence of oriented intervals $\{I_i\}$ on $\alpha_0.$ Choose flute surfaces $F_i$ so that for $j>0,\, |\alpha_{i,j}|=R,$ the value defined at the beginning of Section \ref{R} and so that all finite ends, except for $A^-,$ are punctures. Then the quilted surface $S_{\infty}=S(\alpha_0, p, C,\{F_i\})$ satisfies the hypotheses of the previous section and therefore, Proposition \ref{onepointcritical} holds. 

Let $\{S_n\}$ be the collection  of surfaces converging to $S_{\infty} $  and let  $\Gamma_0\subset\Gamma_1\subset\ldots\subset \Gamma_{\infty} $ be the associated sequence of Fuchsian groups where $\Bbb{H}^2/\Gamma_n=S_n$ as in Lemma \ref{convergeL}.  It follows from  Lemma \ref{convergeL}   and  Proposition \ref{embed}, that for $0\leq k\leq \infty,$  the covering maps $\pi_k:\HH^2\rightarrow \HH^2/\Gamma_k=S_k$ maps $I$ to $\alpha_0$, $\tilde{C}$ to $C$ and $\tilde{p}=ie^{\frac{a}{2}}$ to $p.$

Let $\Gamma_{\infty}^{\varphi}$ be the Fuchsian group $\varphi^{-1}\Gamma_{\infty}\varphi$ and let $\tilde{p}^{\varphi}=\varphi^{-1}(\tilde{p}).$ Recall that $D(\tilde{p} , \Gamma_{\infty} )$ is the Dirichlet polygon of $\Gamma_{\infty}$ centered at   $\tilde{p}$  and 
$\partial_{\infty}D(\tilde{p} , \Gamma_{\infty} )$ is its boundary at infinity.

Recall from Proposition 4 that  one of $\delta_c^s$ or $\delta_c^l$ is critical.  Let $\delta_c$ be one of these rays that is critical.

 Theorem \ref{partial} is an immediate consequence of the following theorem.

\begin{thm}\label{4}
$\partial_{\infty}D(\tilde{p}^{\varphi}, \Gamma_{\infty}^{\varphi})$ consists of the union of the set $C^*$, the interval $\varphi^{-1}([-e^a, -1])$ and a countable set of isolated parabolic fixed points of $\Gamma_{\infty}^{\varphi}$.
\end{thm}

\begin{proof}
Since    $\varphi^{-1}(\partial_{\infty}D(\tilde{p} , \Gamma_{\infty} ))=  \partial_{\infty}D(\tilde{p}^{\varphi}, \Gamma_{\infty}^{\varphi}),$   it   suffices to prove that $\partial_{\infty}D(\tilde{p} , \Gamma_{\infty} )$ consists of the union of the set $\varphi(C^*)$, the interval $ [-e^a, -1] $ and a countable set of isolated parabolic fixed points of $\Gamma_{\infty}.$  
So without loss of generality we suppose $\varphi(z)=z.$

Since $\tilde{p}$ is chosen to be the hyperbolic midpoint of the arc of $I$ with endpoints $i$ and $ie^a$, the Dirichlet polygon $D(\tilde{p},\Gamma_0)$ is bounded by the two geodesics $\psi(1)$ and $\psi(e^a).$ Thus $\partial_{\infty}D(\tilde{p},\Gamma_0)=[-e^a,-1]\cup [1,e^a].$
It follows from Lemma \ref{convergeL} that for each of the subgroups $\Gamma_k,\,   k\in \NN\cup\{\infty\}$ the left half-plane is precisely invariant under the subgroup $\Gamma_0\subset\Gamma_k$; that is, $g(H^-)=H^-$, for $g\in \Gamma_0$ and $g(H^-)\cap H^-=\emptyset$ for $g\in \Gamma_k\setminus\Gamma_0.$ Thus the boundary at infinity of   $D(\tilde{p},\Gamma_k)\cap H^-$ is $[-e^a,-1]$, for all $k\in \NN\cup\{\infty\}.$  

As observed earlier, there is an intimate relationship between points on the boundary at infinity of the Dirichlet polygon at $\tilde{p}$ for $\Gamma_{\infty}$ and the critical rays on $S_{\infty}$ beginning at $p$. In particular, $q\in \partial_{\infty}D(\tilde{p} , \Gamma_{\infty} )$ if and only if
there is a critical ray $\alpha$ on $S_{\infty}$ beginning at $p$ and a lift $\tilde{\alpha}$ to  $
\HH^2$ beginning at $\tilde{p}$ so that $\d{\lim_{t \rightarrow\infty}\tilde{\alpha}(t)=q.}$ In order to simplify notation we shall refer to $q$ as the endpoint of the ray $\tilde{\alpha}.$

For example, the points of $[-e^a,-1]\subset  \partial_{\infty}D(\tilde{p} , \Gamma_{\infty} )$ are the endpoints of lifts of critical rays  beginning at $p$ that go out the finite end on the subsurface $A^-$.

 For each $c\in C,\, \delta_c$ lifts to a geodesic ray $\tilde{\delta}_c$ beginning at $\tilde{p}.$ Since $\delta_c$ is a critical ray, 
the endpoint $q$ of $\tilde{\delta}_c$ lies in the boundary of $ D(\tilde{p} , \Gamma_{\infty} )$.  Let $Q$ denote the set of all endpoints of the lifts of  rays $\tilde{\delta}_c$ for $c\in C.$ We show that $Q=C^*.$ 

Observe that the region $\Omega\subset S_{\infty}$ bounded by  the geodesic rays $\delta_c$ and $\overline{\delta}_c$ is simply connected. Chose the lift $\tilde{\Omega}$ of $\Omega$ to $\HH^2$ so that $\delta_c$ lifts to  $\tilde{\delta}_c.$ Then the boundary of  $\tilde{\Omega}$, $\partial \tilde{\Omega},$ contains a lift of the arc of $\alpha_0\subset \overline{\delta}_c$ passing through $\tilde{p}$ whose length is less than or equal to $\frac{a}{2}$. This lift must then be an arc of $I$ lying between $i$ and $ie^a$. Similarly, $\sigma_c\subset \partial\Omega$, lifts to a geodesic ray $\tilde{\sigma}_c\subset \partial \tilde{\Omega},$ which is orthogonal to the lift $I$ of $\alpha_0.$   The geodesic rays  $\tilde{\delta}_c$ and  $\tilde{\sigma}_c$ are asymptotic and therefore share the same endpoint $q.$  It follows that $\tilde{\sigma}_c(0)=\tilde{c}\in \tilde{C}, $ and then $q=\Lambda^{-1}(\tilde{c})\in C^*.$

Conversely, it follows by the same considerations that, if $c^*\in C^*$ and $c=\pi\circ\Lambda(c^*)$, then the lift $\tilde{\delta}_c$
of $\delta_c$ beginning at $\tilde{p}$ has endpoint $c^*.$ Thus $C^*=Q.$

We next show that if $q\in \partial_{\infty}D(\tilde{p} , \Gamma_{\infty} )$ is not in $[-e^a,-1]$ or $C^*$ then $q$ is an isolated parabolic fixed point of $\Gamma_{\infty}.$ Suppose $\tilde{\alpha}:[0,\infty)\rightarrow \HH^2$ is  a geodesic ray with initial point $\tilde{p}$ and endpoint $q $. Then  $\tilde{\alpha}$ projects to a critical ray $\alpha$ on $S_{\infty}.$ If $\alpha$ goes out the infinite end of $S_{\infty}$, then by Proposition \ref{onepointcritical}, $\alpha=\delta_c$ for some $c\in C$ and then, by the earlier arguments, $q\in C^*.$ Furthermore, if $\alpha$ goes out the end of $S_{\infty}$ corresponding to the annular region $A^-,$ then we have seen that $q\in [-e^a,-1].$ The only remaining possibility is that $\alpha$ goes out a finite end $E$  of $S_{\infty}$, corresponding to a puncture. In that case $q$ must be a parabolic fixed point of $\Gamma_{\infty}.$

It is well known that if $q$ is a parabolic fixed point of $\Gamma_{\infty}$ in $\RR,$  then there is a horocycle (open disc) $N$ in $\HH^2$, tangent to $\RR$ at $q$ so that any geodesic ray in $\HH^2$ that intersects $N$ projects to a self-intersecting geodesic on $S_{\infty}$. Since critical rays are simple, no lift of a critical ray on $S_{\infty}$ beginning at $p$ to one beginning at $\tilde{p}$ can intersect $N$. Consequently, there is a neighborhood of $q$ in $\RR$ that does not contain other boundary points of $D(\tilde{p} , \Gamma_{\infty} ).$

Finally, since $\Gamma_{\infty}$ is countable, there can be only countably many parabolic fixed points. That completes the proof. 
\end{proof}

\subsection{The proof of Theorem \ref{C}}

 Choose          $a<1$ so that $|\alpha_0|=a$ and
    $\sinh^{-1}  [(
\sinh \frac{|\alpha_0|}{2})^{-1}  ] >a.$     
Then the collar neighborhood $C_{S_{\infty}}(\alpha_0, a)$ of $\alpha_0$ of width $2a$ is an embedded annulus in $S_{\infty}=S(\alpha_0, p, C,\{F_i\}).$  There is a unique hyperbolic, twice punctured disc $D$, for which the boundary geodesic $\beta$ has length $a$. As above, the Collar Lemma also tells us that $\beta$ has an embedded collar neighborhood of width $2a$ in $D$.

On $S_{\infty}$ remove the annular region $A^-$ to the left of $\alpha_0$, to get the surface $\overline{S}_{\infty}$ with geodesic boundary $\alpha_0.$ Similarly, on $D$ remove the annular region in the complement of $\beta$ on $D$ to get a twice punctured disc $\overline{D}$ with geodesic boundary $\beta.$ Choose a point $q$ on $\beta$ so that the unique minimal length geodesic arc on $\overline{D}$, separating the punctures, with both endpoints on $\beta$ begins at $q$.

We can now glue the surfaces $\overline{S}_{\infty}$ and $\overline{D}$ together by identifying $\beta$ and $\alpha_0$ so that $p$ and $q$ are matched by the identification. The resulting surface $M_{\infty}$ is a hyperbolic surface with infinitely many ends: one   is an infinite end and the rest are punctures. Let $G$ be the Fuchsian group representing  $M_{\infty}$, which may be chosen so that $G\supset \Gamma_{\infty}.$

The theorem will be proved by showing that $\partial_{\infty}D(\tilde{p},G)$ consists of the set $\varphi(K),$ where $\varphi$ is as in the proof of Theorem \ref{4}, and a set of isolated parabolic fixed points of $G$. Reverting to familiar notation, let $C^*=K $  and, without loss of generality, we suppose that $\varphi(z)=z.$  We consider the problem intrinsically on  $M_{\infty}$. Since $\overline{S}_{\infty} $ sits naturally inside  
$M_{\infty}$, the geodesic rays $\delta_c$  are all  defined on $M_{\infty}$. We shall extend Proposition \ref{onepointcritical} to apply to the surfaces $M_{\infty}$, by showing that the only infinite critical rays on $M_{\infty}$ are still the $\delta_c$ for $c\in C$.

Let us see how the above will suffice to prove the theorem. First, observe that on the new surface there is no finite end of the second kind. As a result, every point in $\partial_{\infty}D(\tilde{p},G)$ is the endpoint of a lift of an infinite critical ray or, as we have seen from earlier arguments, an isolated parabolic fixed point of $G$.
As a consequence of the extended version of Proposition \ref{onepointcritical}, the only lifts of infinite critical rays beginning at $\tilde{p}$ are the $\tilde{\delta}_c$ for $c\in C$, and their endpoints comprise exactly the set $C^*.$

It remains for us to prove, what we call, the extended proposition.  This is done by showing that if $\alpha$ is a critical ray on $M_{\infty},$ then $\alpha\subset \overline{S}_{\infty} \subset M_{\infty}.$ But if $\alpha\subset \overline{S}_{\infty}\subset S_{\infty}$ is an infinite critical ray, then by Proposition \ref{onepointcritical}, $\alpha$ is one of the rays $\delta_c,\,c\in C.$ Suppose $\alpha:[0,u^*]\rightarrow M_{\infty}$ is a geodesic ray with $\alpha(0)=p$, $\alpha\not\subset \overline{S}_{\infty}$ and $\alpha$ goes out the infinite end of $M_{\infty}$. Then $\alpha$ will have non-empty intersection with the interior of the surface $\overline{D}\subset M_{\infty}$ but must eventually lie on $\overline{S}_{\infty}.$ In order for this to occur, there would need to be   values $u_1, u_2, 0\leq u_1  <u_2$ so that  $\alpha([u_1,u_2])\subset \overline{D}$ and $\alpha(u_1), \alpha(u_2)\in \alpha_0.$

 It is now possible to construct a  piecewise geodesic $\overline{\alpha}$ lying entirely on $\overline{S}_{\infty},$ which is strictly shorter than $\alpha$ between the same endpoints.  The geodesic $\overline{\alpha}$ is made by following an arc of $\alpha$ from $\alpha(0)$ to $\alpha(u_1)$, then following an arc of $\alpha_0$ from $\alpha(u_1)$ to $\alpha(u_2)$ and finally, following the arc of $\alpha$ from $\alpha(u_2)$ to $\alpha(u^*).$
 Since the arc $\alpha([u_1,u_2])$ crosses half the collar $C_{M_{\infty}}(\alpha_0, a)$ twice in the interior of $\overline{D}$, its length must be greater than $2a$. But the arc of $\alpha_0$ replacing it has length
less than $a$.   
 Thus $\alpha$ cannot be critical. That completes the proof of Theorem \ref{C}.

%Sigma%%%%%%%%%%%%%%%%%%%%%%%%%%%%%%%%%%%%%%%%

\section{Critical and subcritical rays   get close to $\Sigma$}\label{close}

Consider the set $\Sigma=\{\sigma_c|\ c\in C\}$. Let $\epsilon>0
$ be given and let $N(\Sigma,\epsilon)$ be the
$\epsilon$-neighborhood of $\Sigma.$ We are interested in showing
that the infinite critical and subcritical rays eventually lie in an
$\epsilon$-neighborhood of $\Sigma.$ 

\begin{thm}\label{Sigma}
Let $\sigma:[0,\infty)\rightarrow S_{\infty}$ be an infinite critical or
subcritical ray on $S_{\infty}$. Given $\epsilon>0$ there is a
positive real $t_{\epsilon}$ so that if  $t>t_{\epsilon},$ then $
\sigma(t)\in N(\Sigma,\epsilon).$
\end{thm}

Let $\sigma : [0,\infty)
\longrightarrow S_\infty$ be a geodesic ray that goes out the infinite end of $S_{\infty}.$ We employ a
construction that will enable us to examine a sequence,
$\{\lambda_i\},$ of points on $\sigma$ that lie as
far as possible from $\Sigma,$ in particular, farther than a given positive constant $\epsilon.$  The proof of Theorem \ref{Sigma} will amount to showing that the hypothesis that $\{\lambda_i\}$ is an infinite sequence leads to the conclusion that $\sigma$ is neither  critical nor subcritical, 

  Let $\epsilon >0$ be given.   Define $f:[0,\infty) \longrightarrow [0,\infty)$ by
$f(\lambda)= d(\sigma (\lambda),\Sigma).$ 
%which, for fixed $\lambda $, is the infimum of $ d(\sigma(\lambda),s)$  for $s\in\Sigma.$ 
Clearly, $f$ is continuous, so the set $E=f^{-1}((0,\infty))$ is a countable union of disjoint, open, connected components.  
%Consecutive components of $E$ in $[0,\infty]$ are separated by points $\lambda$ where $\sigma(\lambda)$ crosses $\Sigma.$  
Let $\{V_i\ |\ i \in I^\prime\}$, where $I^\prime$ is some countable index set, be the collection of these open intervals, and let $I \subset I'$ be the set of indices for which $V_i \cap f^{-1}((\epsilon, \infty))$ is nonempty. 
Let $I_b\subset I$ denote the set of $i\in I$ for which $V_i$ is bounded.
For each $i \in I_b,$ let $d_i=max_{\lambda \in \overline{V}_i} f(\lambda)$, and choose exactly one point 
$\lambda_i \in f^{-1}(d_i)\bigcap V_i.$  Clearly, for each $i \in I_b,$  $\sigma(\lambda_i)$ is a choice of a point on $\sigma$ that is farthest from $\Sigma$ for all points $\sigma(\lambda)$, $\lambda \in V_i.$

Given $\lambda, \mu \in L:=\{\lambda_i : i \in I_b\}$, note that if the set $
[\lambda, \mu ] \cap L$ has at least $m$ elements, then
 $\mu-\lambda > 2(m-1)\epsilon.$  As $\mu-\lambda$ is the length of the segment of $\sigma$ that goes from $\sigma(\lambda)$ to $\sigma(\mu)$,  it follows
that for all $\lambda, \mu \in L$, $\lambda < \mu$, $ [\lambda,
\mu ] \cap L$ is a finite set.

Assume for the remainder of this discussion that $L$ is infinite.  Given the remark above, we may now (re-)order $L$ into a strictly increasing sequence
$\{\lambda_i\}_{i=1}^\infty$ where $\lambda_i \rightarrow \infty$ as $i \rightarrow \infty.$  Evidently, the hypothesis that $L$ is an infinite set
implies that each of the (reordered) intervals $V_i$ has finite length and therefore, that $I_b=I.$

For each $i \in \Bbb{N}$ let $\tau_{2i-1}$ and
$\tau_{2i}$ be, respectively, the left and right endpoints of the
interval $\overline{V}_i.$  Note that for all $i,$
$\sigma(\tau_{2i-1})$ and  $\sigma(\tau_{2i})$ lie on scaffolding geodesics in $\Sigma.$  Thus, for each $i \in \Bbb{N},$ there are, not necessarily distinct, scaffolding curves $\beta_{2i-1}$ and $\beta_{2i},$ as well as positive reals $t_{2i-1}$ and $t_{2i},$  such that $\beta_{2i-1}(t_{2i-1})= \sigma(\tau_{2i-1})$ and $\beta_{2i}(t_{2i})= \sigma(\tau_{2i})$. Observe that $\beta_{2i-1}(t_{2i-1})= \sigma(\tau_{2i-1})$ and $\beta_{2i}(t_{2i})= \sigma(\tau_{2i})$ are the endpoints of the smallest connected segment of $\sigma$ that contains $\sigma(\lambda_i),$ and has endpoints lying in $\Sigma.$  In other words, the segment $\sigma((\tau_{2i-1},\tau_{2i}))$ of $\sigma,$ which includes the point $\sigma(\lambda_i),$ has empty intersection with $\Sigma.$ Since each point of $S_{\infty}$ lies either on $\Sigma$  or in the interior of a cut-open subflute, it follows that $\sigma((\tau_{2i-1},\tau_{2i}))$ lies on a cut-open subflute, denoted $F^*_i.$ Now, $\beta_{2i-1}$ is one of the boundary scaffolding geodesics of $F^*_i;$  let $\beta'_i$ be the other and observe that either $\beta'_i=\beta_{2i-1}$ or $\beta'_i=\beta_{2i}.$   
Note that, from the construction, for each $i \in \Bbb{N},$ we have $\sigma:[\tau_{2i-1},\tau_{2i}]\rightarrow F^*_i,$ $\sigma(\tau_{2i-1}),\,\sigma(\tau_{2i})\in\, \beta_{2i-1}\cup\beta'_i=\partial F^*_i,$ and $\lambda_i \in [\tau_{2i-1},\tau_{2i}]$ such that $d(\sigma(\lambda_i),(\beta_{2i-1} \cup \beta'_i)) \ge d(\sigma(\lambda),(\beta_{2i-1} \cup \beta'_i))$ for all $\lambda \in [\tau_{2i-1},\tau_{2i}],$  that is, $\sigma(\lambda_i)$ realizes the maximal value for $d(\sigma(\lambda), \beta_{2i-1}\cup\beta'_i)$  among all points on   $\sigma([\tau_{2i-1},\tau_{2i}])$.  By Lemma \ref{newFive}, $\sigma(\lambda_i)$ must lie on a $\gamma$-curve of $F^*_i.$  Applying Lemma \ref{newSix} to this arc of $\sigma$ on the flute $F^*_i$, we obtain
$$\tau_{2i}-\tau_{2i-1}>|t_{2i}-t_{2i-1}|+2 \log\left (\frac {k^2+1}{2k}\right ),$$
where $\log k=d(\sigma(\lambda_i), \beta_{2i-1}\cup\beta'_i).$
 
As earlier, let $\kappa(R)=2 \log\big (\frac {k^2+1}{2k}\big ),$ where $\log k = R.$ For each $i\in\NN$ set $\epsilon_i=d(\sigma(\lambda_i), \beta_{2i-1}\cup\beta'_i).$ Then from the construction we have $\epsilon_i\geq\epsilon$ for all $i\in\NN.$ Since $\kappa$ is strictly increasing for $R\geq 0,$ we have proved the first part of the next lemma; the second part is a direct consequence of Lemma \ref{forSigma}.

\begin{lemma}\label{cepsilon}
Let $\epsilon>0$ be given.  With definitions as above,
 there is a positive constant $\kappa(\epsilon)$ such that
$\tau_{2i}-\tau_{2i-1}>
|\ t_{2i}-t_{2i-1} \ |+\kappa(\epsilon).$  Generally, for all $i$,
$\tau_{i+1}-\tau_i > |\ t_{i+1}-t_{i} \ |$; in particular,
$\tau_{2i+1}-\tau_{2i} > |\ t_{2i+1}-t_{2i} \ |.$
\end{lemma}

\noindent {\em Proof of Theorem \ref{Sigma}. }\, 

Let $\epsilon >0$ be given.  For the geodesic ray $\sigma,$ construct the collection $V=\{V_i\ |\ i \in I\}$ as defined above. Recall that each for each $V_i,$ $\sigma(V_i)$ lies in the cut-open, untwisted flute $F^*_i$ and $V_i$ contains a point $\lambda_i$ for which $d(\sigma(\lambda_i),\Sigma)>\epsilon.$ Further, any point $t>0$ for which $d(\sigma(t),\Sigma)>\epsilon$ must lie in some $V_i.$   Therefore, if the collection $V$ is a finite collection of sets, and each $V_i$ in the collection is bounded, then the conclusion of the theorem is true, namely, there is a $t_\epsilon>0$ such that, $\sigma(t)\in N(\Sigma,\epsilon)$ for $t> t_\epsilon.$   

We note first that $V_i$ cannot be unbounded.  In that case, since $\sigma(V_i)$ lies in exactly one cut-open, untwisted flute, $F^*_i,$ the curve $\sigma$ would, after a point, lie entirely in the flute $F^*_i.$  But then by the hypothesis that the each  flute $F_i$ is of the first kind, the results in \cite{haas} can be used to deduce that $\sigma$ must eventually lie arbitrarily close to the scaffolding curves bounding this subflute, and therefore must eventually lie in $N(\Sigma,\epsilon). $

We now know the sets $V_i$ in the collection $V$ are each bounded.  It remains to show that the collection is finite.  Suppose that $V$ is an infinite collection.  We will show that in this case, $\sigma$ is neither critical not subcritical.  Employing the earlier construction, there are infinite sequences of scaffolding curves $\beta_i,$ $\beta'_i$ and infinite sequences
$\lambda_i$, $\tau_i$, $t_i$ $\in [0,\infty],$ $i \in \Bbb{N}$ as above, where for each $i,$ $d(\Sigma,
\sigma(\lambda_i))>\epsilon$; in particular,
$d(\sigma(\lambda_i), (\beta_{2i-1}\cup \beta'_{i}))>\epsilon$. By
Lemma \ref{cepsilon},  there is a positive constant $\kappa(\epsilon)$
such that, for each $i \in \Bbb{N}$, $$\tau_{2i}-\tau_{2i-1} > |\
t_{2i}-t_{2i-1} \ |+\kappa(\epsilon);$$ and for all $i,$
$$\tau_{2i+1}-\tau_{2i} > t_{2i+1}-t_{2i}.$$

Choose $N \in \Bbb{N}$ such that $N\kappa(\epsilon) >
|\alpha_0|+d(\sigma(0),\alpha_0)+\ |t_1-\tau_1\ | +m,$ where
$m$ is an arbitrarily chosen positive real.
 The length of the curve
$\sigma([0,\tau_{2N}])$  is $\tau_{2N}$.
 From Lemma \ref{cepsilon} we have,

\begin{align}
\tau_{2N}-\tau_1
&=\Sigma_{i=1}^{2N-1}(\tau_{i+1}-\tau_{i})\nonumber\\
&=\Sigma_{i=1}^{N}(\tau_{2i}-\tau_{2i-1})+
\Sigma_{i=1}^{N-1}(\tau_{2i+1}-\tau_{2i})\nonumber\\
&>\Sigma_{i=1}^{N}(|\ t_{2i}-t_{2i-1}\ | +\kappa(\epsilon)) + \Sigma_{i=1}^{N-1}|\
t_{2i+1}-t_{2i}\ |\nonumber\\
&\geq N\kappa(\epsilon) + \Sigma_{i=1}^{N}(
t_{2i}-t_{2i-1}) + \Sigma_{i=1}^{N-1}(t_{2i+1}-t_{2i})\nonumber\\
&=N\kappa(\epsilon)+ \Sigma_{i=1}^{2N-1}(t_{i+1}-t_i)\nonumber\\
&=N\kappa(\epsilon)+ t_{2N}-t_1\nonumber\\
&> |\alpha_0|+d(\sigma(0),\alpha_0)+\ |t_1-\tau_1|+m +t_{2N}-t_1.\nonumber  
\end{align}

Therefore,
\begin{align}
\tau_{2N}&> |\alpha_0|+d(\sigma(0),\alpha_0)+|t_1-\tau_1|-(t_1-\tau_1)+t_{2N}+m \nonumber\\
&\geq |\alpha_0|+d(\sigma(0),\alpha_0)+t_{2N}+m.\nonumber  
\end{align}

Note that the piecewise curve beginning at $\sigma(0)$, proceeding along a minimal length geodesic  to a point on $\alpha_0$, then along $\alpha_0$ to the scaffolding curve $\beta_{2N}(0)$, along $\beta_{2N}$ to
$\beta_{2N}(t_{2N}),$ has length less than $|\alpha_0|+d(\sigma(0),\alpha_0)+t_{2N}$.  It follows from the inequality above that the length of the curve $\sigma(0,\tau_{2N})$ is
larger than this piecewise curve by at least $m,$ where $m$ was arbitrary.   
Therefore, it follows that $\sigma$ can be neither critical nor subcritical.
$\hfill\Box$\\
\\

%End Sigma %%%%%%%%%%%%%%%%%%%%%%%% 

\end{document}